\documentclass[10pt,twoside]{amsart}

\usepackage{geometry}
\usepackage{pdfsync}
\usepackage{hyperref}
\usepackage{amsmath,amssymb,esint,amsthm,bbm}
\usepackage[active]{srcltx}

\usepackage[dvips]{graphicx}
\catcode`@=11 \@addtoreset{equation}{section} \catcode`@=12

\title[Existence and Spectral theory]
{Existence and Spectral theory for Weak Solutions of Neumann and Dirichlet problems for linear
degenerate elliptic operators with rough coefficients}

\email{dario.monticelli@gmail.com (Dario Daniele Monticelli)}
\email{scott.rodney@gmail.com (Scott Rodney)}


\subjclass[2010]{}

\newcommand{\bu}{{\bf u}}
\newcommand{\bv}{{\bf v}}
\newcommand{\e}{\varepsilon}
\newcommand{\gb}{\mathfrak{b}}
\newcommand{\vphi}{\varphi}

\newcommand{\field}[1]{\mathbb{#1}}
\newcommand{\R}{\field{R}}

\newcommand{\N}{\field{N}}

\newcommand{\crochet}[2]{\left\langle{#1},\,{#2}\right\rangle}

\newcommand{\wt}{\widetilde}

\newcommand{\bea}{\begin{eqnarray}}
\newcommand{\cend} {\end{center}}

\newcommand{\eea}{\end{eqnarray}}

\newcommand{\ra}{\rightarrow}

\newcommand{\beginc} {\begin{center}}

\newtheorem{thm}{\textbf{Theorem}}[section]
\newtheorem{lem}[thm]{\textbf{Lemma}}
\newtheorem{pro}[thm]{\textbf{Proposition}}
\newtheorem{rem}[thm]{\textbf{Remark}}
\newtheorem{cor}[thm]{\textbf{Corollary}}
\theoremstyle{remark}

\newtheorem{exe}[thm]{\textbf{Example}}
\theoremstyle{definition}
\newtheorem{defn}[thm]{\textbf{Definition}}
\newtheoremstyle{Claim}{}{}{\itshape}{}{\itshape\bfseries}{:}{ }{#1}
\theoremstyle{Claim}

\begin{document}

\maketitle \vspace{-0,3cm}

\begin{center}
\textsc{D. D. Monticelli\footnote{Universit\`{a} degli Studi di
Milano, Italy. Partially supported by GNAMPA project with title ``Equazioni differenziali con invarianze in analisi globale'', by GNAMPA section ``Equazioni differenziali e sistemi dinamici'' and by MIUR project ``Metodi variazionali e topologici nello studio di fenomeni nonlineari''} and S. Rodney\footnote{Cape Breton University, Sydney, NS Canada.  Partially supported by The Natural Sciences and Engineering Research Council of Canada (NSERC) Discovery Grants program.}}
\end{center}


\begin{abstract}
In this paper we study existence and spectral properties for weak
solutions of Neumann and Dirichlet problems associated to second
order linear degenerate elliptic partial differential operators $X$,
with rough coefficients of the form $$X=-\text{div}(P\nabla )+{\bf
HR}+{\bf S^\prime G} +F$$ in a geometric homogeneous space setting
where the $n\times n$ matrix function $P=P(x)$ is allowed to
degenerate.  We give a maximum principle for weak solutions of
$Xu\leq 0$ and follow this with a result describing a relationship
between compact projection of the degenerate Sobolev space
$QH^{1,p}$ into $L^q$ and a Poincar\'e inequality with gain adapted
to $Q$.
\end{abstract}

\section{Introduction}

This paper studies existence of weak solutions and spectral
properties for Dirichlet and Neumann problems on bounded domains of
$\R^n$ for linear second order degenerate elliptic partial
differential operators with rough measurable coefficients of the
form \bea\label{I0.1} X &=& -\text{div}(P\nabla )+ {\bf HR} + {\bf
S'G} + {\bf F}, \eea where $P=P(x)$ is a nonnegative definite
symmetric measurable matrix function, ${\bf H,G,F}$ are vector
valued functions and ${\bf R,S}$ are collections of first order
subunit vector fields. We refer the reader to section $1.6$ for a
precise description of the constituents of $X$.

Before continuing we briefly describe how the paper is organized.
Our results are developed in a general axiomatic framework similar
to those used in \cite[Section 3]{CRW}, \cite{SW1,SW2},
\cite{MRW}, that we outlined in Sections 1.1 through 1.6. This
axiomatic setting includes the definition of geometric homogeneous
space, of degenerate Sobolev spaces associated to a nonnegative
definite symmetric measurable matrix function comparable to
$P$ and also gives Poincar\'e--Sobolev type inequalities. Existence of weak
solutions and spectral properties associated to the Neumann problem
for $X$ (referred to as the $X$-Neumann problem) are studied in
Section $2$ and the main results are Theorems \ref{thm2},
\ref{thm3}, Corollary \ref{cor3}, and Theorem \ref{thm6}.  We also
mention that we develop a helpful example (see Example \ref{exe1})
in full detail for the reader's convenience.  We give a spectral
result for the Dirichlet problem for $X$ (referred to as the
$X$-Dirichlet problem) in Section $3$, see Theorem \ref{thm7}, that
compliments the existence results of \cite{ScottCJM}.  Section 4
contains a maximum principle for weak solutions of $Xu\leq 0$ and in Section $5$ we demonstrate a relationship between compact embeddings of Sobolev spaces and global Poincar\'e inequalities with
gain; we refer the reader to Theorems \ref{thm8} and \ref{thmP} for
these results.  All of our results are developed in the spirit of
\cite{Evans} and \cite{GT} using ideas presented in \cite{CRW},
\cite{DarioPhD}, \cite{DarioH=W}, \cite{DarioJEMS}, \cite{MRW},
\cite{ScottPhD}, \cite{Scott}, \cite{ScottCJM}, \cite{SW1},
\cite{SW2}, and other related works.


\subsection{Geometric Homogeneous Spaces:} Let $\Omega$ be a bounded open subset of $\mathbb{R}^n$ and $\rho$ a symmetric quasi-metric defined in $\Omega$.  More precisely, $\rho:\Omega\times\Omega\ra [0,\infty)$ and there is a $\kappa\geq 1$ so that for every $x,y,z\in \Omega$ each of the following are satisfied:
\bea\label{I1}\bullet && \rho(x,y)=0 \iff x=y,\nonumber\\
\bullet && \rho(x,y) = \rho(y,x),\text{ and} \nonumber\\
\bullet && \rho(x,y) \leq \kappa[\rho(x,z) + \rho(z,y)].
\eea
Given $x\in\Omega$ and $r>0$ we define the sets
\bea B(x,r) &=& \{y\in\Omega\;:\;\rho(x,y)<r\},\text{ and}\nonumber\\
D(x,r) &=&\{y\in\Omega\;:\; |x-y|<r\}.\nonumber
\eea
We always refer to $B(x,r)$ as the quasimetric ball centered at $x$ with radius $r$.  $D(x,r)$ is the corresponding Euclidean ball with the same center and radius. \begin{defn} We say that the collection of quasimetric balls ${\mathcal{B}}=\{B(x,r)\}_{r>0;x\in\Omega}$ is \emph{locally geometrically doubling} if given any compact $K\subset\Omega$ there is a $\delta=\delta(K)>0$ so that for all $x\in K$ and $0<r'\leq r<\delta$, the $\rho$-ball $B(x,r)$ may contain at most $C=C(r/r')$ disjoint $\rho$-balls of radius $r'$ where $C$ is independent of both $x$ and $K$.
\end{defn}
It is important to note that the notion of local geometric doubling
is weaker than local doubling for Lebesgue measure on the collection
of quasimetric balls centered in $\Omega$.  Recall that a
measure $\mu$ is locally doubling for the collection $\mathcal{B}$
if given a compact set $K\subset \Omega$ there is a
$\delta_1=\delta_1(K)$ so that for any $x\in K$ and $0<r<\delta_1$
we have $|B(x,2r)|_\mu \leq C|B(x,r)|_\mu$ where $C>0$ is
independent of $x,r,K$.  Helpful discussions on geometric doubling
and doubling conditions are found in \cite{HyK}, and \cite{HyM}. We
now define the topological space on which we build our results.
\begin{defn} The pair $(\Omega,\rho)$ is called a \emph{geometric homogeneous space} if $\rho$ is a symmetric quasimetric defined in $\Omega$ for which the three following properties hold.
\bea\label{Cond1} \bullet && \text{For all $x\in\Omega$ and $r>0$ the quasimetric $\rho$-ball $B(x,r)$ is an open set}\nonumber\\
&&\text{with respect to Euclidean topology.}\nonumber\\
\bullet && \text{For all }x,y\in \Omega, \;|x-y|\ra 0 \text{ if } \rho(x,y)\ra 0.\\
\label{I3}\bullet && \text{The collection of $\rho$-balls }\mathcal{B}=\{B(x,r)\}_{x\in\Omega,r>0}\text{ is locally}\nonumber\\
&&\text{geometrically doubling.}\nonumber
\eea
\end{defn}

\begin{rem} \hspace{1in}
\begin{enumerate}
\item The converse of (\ref{Cond1}) holds automatically since $\rho$-balls are open sets with respect to Euclidean topology.
\item  In many situations it is convenient to work with $\rho$-balls that do not intersect the boundary of $\Omega$ and is the reasoning behind condition (\ref{Cond1}).  To see this notice that (\ref{Cond1}) may be restated as: given $x\in\Omega$ and $\epsilon>0$ there is a positive $\delta''=\delta''(x,\epsilon)$ so that $B(x,\delta) \subset D(x,\epsilon)$.  From this, given $x\in\Omega$ there is a $t_0>0$ so that $\overline{B(x,r)}\subset\Omega$ for all $0<r<t_0$.

\end{enumerate}
\end{rem}

\subsection{Functional Spaces}  Fix an open $\Theta\subset\Omega$, and a non-negative definite $n\times n$ measurable matrix $Q(x)$ with $|Q(x)|\in L^\infty_{loc}(\Omega)$.  The degenerate Sobolev space $\mathcal{W}_Q^{1,2}(\Theta)$ is the collection of pairs $(f,\vec{g})$ obtained via isomorphism from the space $QH^1(\Theta)$ defined as the completion, with respect to the norm
\bea\label{Wnorm} ||u||_{QH^1(\Theta)} &=& \Big(\int_\Theta |u|^2dx
+ \int_\Theta |\sqrt{Q}\nabla u|^2dx\Big)^\frac{1}{2}, \eea of the
collection \bea\label{lipspace} Lip_{Q}(\Theta) = \{\vphi\in
Lip_{loc}(\Theta)\;:\; |\sqrt{Q}\nabla\vphi| \in L^2(\Theta)\} \eea
of locally Lipschitz functions having finite $QH^1(\Theta)$ norm. In
all of our developments we will denote the vector valued function
$\vec{g}$ of the pair $(f,\vec{g})\in QH^1(\Theta)$ by writing
$\vec{g}=\nabla f$, and we will refer to it as the gradient part (or
simply the gradient) of $f$. The canonical projection map
$i:QH^1(\Theta)\rightarrow L^2(\Theta)$ defined by
$$i\big((f,\vec{g})\big)=f$$ is obviously continuous, but it need not be injective. It should be kept in mind at all times
that $\vec{g}$ need not be uniquely determined by $f$; see
\cite{FKS} for a well known example. In order to keep notation
simple we will often abuse notation by writing $f\in QH^1(\Theta)$
in place of $(f,\nabla f)\in QH^1(\Theta)$. In cases where confusion
may arise, specifically in Section 4, we will use bold faced
characters when referring to elements of $QH^1(\Theta)$; i.e. we
will write ${\bf f}\in QH^1(\Theta)$ in place of $f\in
QH^1(\Theta)$.

The space $QH^1_0(\Theta)\; \big(\text{respectively }\mathcal{W}_{Q,0}^{1,2}(\Theta)\big)$ is obtained in a similar manner, but in this case we complete
the set $Lip_0(\Theta)$, the set of those Lipschitz functions having compact support in $\Theta$, with respect to the norm (\ref{Wnorm}).  Notice that all such functions have finite $QH^1(\Theta)$ norm as $|Q|\in L^\infty_{loc}(\Omega)$. \\
\indent For clarity, we will always write $QH^1(\Theta)$ and
$QH^1_0(\Theta)$ in place of $\mathcal{W}_Q^{1,2}(\Theta)$ and
$\mathcal{W}^{1,2}_{Q,0}(\Theta)$ respectively, taking isomorphism in context. We adopted this notation in lieu of
$W^{1,2}_Q(\Omega)$ and $W^{1,2}_{Q,0}(\Omega)$, as used in
\cite{SW1,SW2,CRW,MRW}, in order to agree with classical literature,
see for example \cite{H=W}, where it is convention that ``$W$"
spaces refer to Sobolev spaces defined with respect to
distributional derivatives.  We also mention that it is possible to
introduce definitions and make similar considerations for the spaces
$QH^{1,p}(\Theta),\;QH^{1,p}_0(\Theta)$ (or $W^{1,p}_Q(\Omega),\;
W^{1,p}_{Q,0}(\Omega)$  as in the above references) for $1\leq p
<\infty$, even in the case where $|Q(x)|$ is locally unbounded.  We
invite the interested reader to see \cite{CRW}, \cite{MRW},
\cite{SW2} for the construction of these spaces and related objects.

\subsection{Sobolev and Poincar\'e Inequalities}

Essential to most of the arguments to follow are Poincar\'e -- Sobolev type inequalities adapted to the matrix $Q$.  We list these inequalities now noting that they are not assumed to hold at all times but, rather, called upon when necessary.  To describe these efficiently, we fix a continuous function $r_1:\Omega\ra (0,\infty)$ to be used as a common radius restriction for quasimetric balls.\\
\indent{\bf The Local Poincar\'e Inequality.}  We say that the local Poicar\'e inequality of order $p$ holds if there are constants $C_2>0$ and $\gb\geq 1$ so that for every $\rho$-ball $B(y,r)$ centered in $\Omega$ with $\gb r\in (0,r_1(y))$ the inequality
\bea\label{Poincare-loc}
\Big(\frac{1}{|B_r|}\int_{B_r} | f-f_{B_r}|^pdx\Big)^\frac{1}{p} &\leq& C_2 r\Big(\frac{1}{|B_{\gb r}|}\int_{B_{\gb r}}|\sqrt{Q}\nabla f |^pdx\Big)^\frac{1}{p}
\eea
holds for all $f\in Lip_{loc}(\Omega)$.  Notice that a continuity argument allows one to extend (\ref{Poincare-loc}) to hold for all pairs $(f,\nabla f)\in QH^{1,p}(\Omega)$.


\textbf{The Global Sobolev Inequality.} For an open set $\Theta\subset\Omega$ with $\overline{\Theta}\subset\Omega$, we say that the global Sobolev inequality holds on $\Theta$ holds if there are positive constants $C_3>0$ and $\sigma>1$ such that
\begin{eqnarray}
\label{GWP}\left(\int_{\Theta}
|f|^{2\sigma}\,dx\right)^{\frac{1}{2\sigma}}&\le& C_3
\left(\int_{\Theta} |\sqrt{Q}\nabla
f|^2\,dx\right)^{\frac{1}{2}}
\end{eqnarray}
holds for all $f\in Lip_0(\Theta)$.  Again, a continuity argument shows that (\ref{sobloc}) also holds for all $(w,\nabla w)\in QH^1_0(\Theta)$ if it holds for all $f\in Lip_0(\Theta)$.
\begin{rem} It is possible to replace (\ref{GWP}) with a local version of the form
\bea  \label{sobloc} \left(\int_{B(x,r)}
|f|^{2\sigma}\,dx\right)^{\frac{1}{2\sigma}}&\le& C'
\left(r^2\int_{B(x,r)} |\sqrt{Q}\nabla f|^2\,dx+
\int_{B(x,r)}|f|^2dx\right)^{\frac{1}{2}} \eea holding for all $f\in
Lip_0(B(x,r))$ and quasimetric balls $B(x,r)$ with $0<r<r_1(x)$.  In
\cite{Scott}, inequality (\ref{GWP}) is proved for any open $\Theta$
such that $\overline{\Theta}\subset\Omega$ under the hypotheses
given in subsections 1.1 and 1.2, provided Lebesgue measure is
doubling for the collection of quasimetric balls $\mathcal{B}$ and
both (\ref{Poincare-loc}) with $p=2$ and (\ref{sobloc}) hold. One may even
replace the doubling assumption in \cite{Scott} with local geometric
doubling; see \cite{CRW,Scott}.
\end{rem}

\textbf{The Global Poincar\'{e} Inequality With Gain $\omega$.} For an open subset $\Theta$ of $\Omega$ satisfying $\overline{\Theta}\subset\Omega$ we say that the global Poincar\'e inequality with gain $\omega>1$ holds on $\Theta$ if there are constants $C_4>0$ and $\omega>1$ such that
\begin{eqnarray}
\label{GP}\left(\int_{\Theta}
|f-f_\Theta|^{2\omega}\,dx\right)^{\frac{1}{2\omega}}&\le& C_4
\left(\int_{\Theta} |\sqrt{Q}\nabla f|^2\,dx\right)^{\frac{1}{2}}
\end{eqnarray}
holds for all $f\in Lip_Q(\Theta)$.

\begin{rem}\label{rem1.5}\hspace{1in}
\begin{enumerate}
\item If the global Poincar\'{e} inequality \eqref{GP} holds, then H\"older's inequality implies that the \emph{{\bf Global Weak Poincar\'e Inequality with gain $\omega>1$}}:
\bea \label{globsob}\left(\int_{\Theta}
|f|^{2\omega}\,dx\right)^{\frac{1}{2\omega}}&\le& C_4
\left(\int_{\Theta} |\sqrt{Q}\nabla f|^2\,dx + \int_\Theta |f|^2\,dx\right)^{\frac{1}{2}}
\end{eqnarray}
also holds for all $f\in Lip_Q(\Theta)$.
\item In the elliptic case ($Q(x)=Id$), inequalities of the form (\ref{GP}) and (\ref{globsob}) are proved when the boundary of $\Theta$ is sufficiently regular.  For example, $\partial\Theta\in C^{0,1}$ is used in \cite{GT} for such purposes.  See \cite[Theorem 7.26]{GT} and related discussions.
\end{enumerate}
\end{rem}
There is a large body of literature discussing the inequalities just mentioned.  We refer the reader to \cite{Evans}, \cite{HK}, and \cite{J} for helpful discussions and examples. \\

As a last remark we mention that it is possible to obtain all the theory to follow if Lebesgue measure is replaced by any Radon measure $\mu$ that is absolutely continuous with respect to Lebesgue measure.  The changes required for this are:
\begin{itemize}
\item replace $|E|$ with $|E|_\mu$ and any almost everywhere considerations shifted to
$\mu-$a.e.;
\item Incorporate $\mu$ into the definition of Sobolev space by replacing $||u||_{QH^1(\Theta)}$ with $$||u||_{QH^1_\mu(\Theta)} = ||u||_{L^2_\mu(\Theta)} + ||\sqrt{Q}\nabla f||_{L^2_\mu(\Theta)};$$
\item replace $dx$ with $d\mu$ in all integrals. 
\end{itemize}

\subsection{Compact Projections for Degenerate Sobolev Spaces.} In our main results, it is essential that we are equipped with a
compact mapping from $QH^{1}(\Theta)$ (respectively
$QH^{1}_0(\Theta)$) into $L^2(\Theta)$ and this serves as our first
application of the inequalities just listed. The following results
are adapted from \cite{CRW} and we refer the reader to that work for
more general statements and weighted results.

\begin{pro}[\cite{CRW} Corollary 3.25]\label{compact1} Let $(\Omega,\rho)$ be a geometric homogeneous space and fix an open set $\Theta$ satisfying $\overline{\Theta}\subset \Omega$.  Suppose that the local Poincar\'e inequality (\ref{Poincare-loc}) with $p=2$ and the global Sobolev inequality (\ref{GWP}) hold.  Then the projection map $i:QH^1_0(\Theta) \ra L^q(\Theta)$ defined by
\bea\label{mapme} i\big( (u,\nabla u)\big) = u \eea is a compact
mapping for all $q\in [1,2\sigma)$.
\end{pro}

\begin{pro}[\cite{CRW} Theorem 3.14]\label{compact2} Let $(\Omega,\rho)$ be a geometric homogeneous space and fix an open set $\Theta$ satisfying $\overline{\Theta}\subset \Omega$.  Suppose that the local Poincar\'e inequality (\ref{Poincare-loc}) with $p=2$ and the global weak Poincar\'e inequality (\ref{globsob}) with gain $\omega>1$ hold.  Then the projection map $i:QH^1(\Theta) \ra L^q(\Theta)$ defined by
\bea\label{mapme1} i\big( (u,\nabla u)\big) = u \eea is a compact
mapping for all $q\in [1,2\omega)$.
\end{pro}

\subsection{Notation.} Consider a vector field
\bea W(x)=\sum_{i=1}^nw_i(x)\frac{\partial}{\partial
x_i}=\big(w_1(x),\ldots,w_n(x)\big)\cdot \nabla\nonumber
\eea
 If $u$ is a real valued function on $\mathbb{R}^n$ and $\nu$
is a vector in $\R^n$ we adopt the notation
$$Wu=\sum_{i=1}^nw_i\frac{\partial u}{\partial
x_i},\qquad\crochet{\nu}{W}=\sum_{i=1}^nw_i\nu_i,$$ where
$\crochet{\cdot\;}{\cdot}$ denotes the standard inner product on
$\R^n$. The formal adjoint $W^\prime(x)$ of the vector field $W(x)$
is defined by
$$W^\prime(x)u:=-\text{div}\big(w_1(x)u(x),\ldots,w_n(x)u(x)\big)=-\sum_{i=1}^n\frac{\partial}{\partial x_i}\big(w_i(x)u(x)\big).$$
A vector field $W(x)$ as above is always identified with the vector valued function\\
$\big(w_1(x),\ldots,w_n(x)\big)$ and is said to be
\emph{subunit} with respect to the matrix $Q$ in $\Theta$ if
\bea \label{subunit} \left(\sum_{i=1}^nw_i(x)\xi_i\right)^2\leq\crochet{\xi}{Q(x)\xi}
\eea
for every $\xi\in\R^n$ and almost every $x\in\Theta$.
\begin{rem}\label{subunitdef} If a vector field $W(x)$ is subunit with respect to the matrix $Q=Q(x)$ in $\Theta$ we will simply refer to it as a ``\emph{subunit vector field}" with the set $\Theta$ and matrix $Q$ taken in context.
\end{rem}

Given $N\in\N$, an $N$--tuple $\mathbf{W}=\big(W_1,\ldots,W_N\big)$
of vector fields and an $\R^N$--valued function
$\mathbf{G}=(g_1,\ldots,g_N)$, $\mathbf{WG}$ denotes the ``inner
product'' of $\mathbf{W}$ and $\mathbf{G}$, i.e.
$$\mathbf{WG}=\sum_{i=1}^NW_i(x)g_i(x).$$
Lastly, If $u$ is a real valued function,
$$\mathbf{GW}u=\sum_{i=1}^Ng_i(x)W_i(x)u(x),\qquad\mathbf{W^\prime}(\mathbf{G}u)=\sum_{i=1}^NW_i^\prime(x)\big(g_i(x)u(x)\big).$$

\subsection{Second Order Linear Degenerate Elliptic Operators with Rough Coefficients}

 Let $\Omega$ and $Q(x)$ be as in Chapter $1$, $\Theta$ be a bounded
domain with $\overline\Theta\subset\Omega$.  We consider operators $X$ defined by equations of the form
\begin{equation}\label{diffop}
Xu=-\text{div}(P\nabla
u)+\mathbf{H}\mathbf{R}u+\mathbf{S}^\prime(\mathbf{G}u)+Fu\text{ in }\Theta,
\end{equation}
where $P=P(x)$ is a bounded measurable nonnegative definite symmetric
$n\times n$ matrix defined in $\overline\Theta$ and comparable with
$Q(x)$ in $\Theta$.  That is, there exist constants $c_1,\,C_1>0$ so that for every $\xi\in\R^n$ and almost every $x\in\Theta$ one has
\begin{equation}\label{4}
c_1\crochet{\xi}{Q(x)\xi}\leq\crochet{\xi}{P(x)\xi}\leq
C_1\crochet{\xi}{Q(x)\xi}.
\end{equation}
$\mathbf{R}$, and $\mathbf{S}$ are, for some $N\in\mathbb{N}$,
$N$--tuples of subunit vector fields with respect to the matrix
$Q(x)$, see \eqref{subunit}. $\mathbf{H}$ and $\mathbf{G}$ are
$\R^N$--valued measurable functions in $\Theta$ and $F$ is a
real-valued measurable function defined in $\Theta$.

\begin{rem}\label{rem0}
The sign convention on the principal part of the operator $X$ is as in \cite{Evans} and \cite{GT}. S. Rodney in
\cite{ScottPhD}, \cite{Scott}, and \cite{ScottCJM} adopted the opposite convention.\
\end{rem}
An operator $X$ as above will always be referred to as a ``second
order linear degenerate elliptic operator with rough coefficients".
We will study existence and spectral properties of both
Dirichlet and Neumann problems associated to $X$ in the context of
two ``negativity" conditions that we now introduce.

\begin{defn}\label{neg1}  We say that a second order linear
degenerate elliptic operator with rough coefficients $X$ satisfies {\bf negativity condition (1)} if and only if either
\begin{itemize}
\item[i)] there exists $\e>0$ such that for every $(u,\nabla u),(v,\nabla v)\in QH^1(\Theta)$ satisfying $uv\geq0$ almost everywhere in
$\Theta$ one has
$$\int_\Theta Fuv+\mathbf{GS}(uv)\,dx\geq\e\int_\Theta
uv\,dx$$
$$\text{or}$$
\item[ii)] for every $(u,\nabla u),(v,\nabla v)\in
QH^1(\Theta)$ satisfying $uv\geq0$ almost everywhere in $\Theta$ one
has
$$\int_\Theta Fuv+\mathbf{HR}(uv)\,dx\geq\e\int_\Theta uv\,dx.$$
\end{itemize}
\end{defn}
\begin{defn}\label{neg2} We say that $X$ satisfies {\bf negativity condition (2)} if and only if either
\begin{itemize}
\item[i)] for every $(u,\nabla u),(v,\nabla v)\in QH_0^1(\Theta)$ such that $uv\geq0$ almost everywhere in
$\Theta$ one has
$$\int_\Theta Fuv+\mathbf{GS}(uv)\,dx\geq0$$
$$\text{or}$$
\item[ii)] for every $(u,\nabla u),(v,\nabla v)\in QH^1_0(\Theta)$ such that $uv\geq0$ almost everywhere in
$\Theta$ one has
$$\int_\Theta Fuv+\mathbf{HR}(uv)\,dx\geq0.$$
\end{itemize}
\end{defn}
It is useful to note that if $X$ satisfies negativity condition (2),
part i), then the operator $X$ is a member of the ``Nonnegative
Class" as described in \cite{ScottCJM}.

\section{Existence and Spectral Results for the $X$-Neumann Problem}
The presentation of the results in this section follows in part the
presentation of the analogous results for second order elliptic
Dirichlet problems on bounded domains of $\R^n$ as given in
\cite{Evans}.  Furthermore, many of the arguments used, vis-a-vis existence of weak solutions, follow similar paths as those found in \cite{ScottCJM}.\\

Let $\Omega$ and $Q=Q(x)$ be as in sections 1.1 and 1.2 respectively.  We begin by describing the action of a subunit vector field on elements and products of elements of $QH^{1}(\Theta)$.



\begin{lem}\label{lem1}
Let $\Theta$ be a bounded open set such that
$\overline\Theta\subset\Omega$, let $(u,\nabla u),\,(v,\nabla v)\in
QH^1(\Theta)$ and let $W(x)=\big(w_1(x),\ldots,w_n(x)\big)$ be a
subunit vector field in $\Theta$. Assume that the global weak
Poincar\'{e} inequality with gain $\omega>1$ holds, see
\eqref{globsob}. Then:
\begin{itemize}
\item[1)] $Wu$ is defined as an element of $L^2(\Theta)$, with
$\|Wu\|_{L^2(\Theta)}\leq\|u\|_{QH^1(\Theta)}$, and
$$W(x)\varphi(x)=\sum_{i=1}^nw_i(x)\frac{\partial\varphi}{\partial
x_i}(x)$$ for every locally Lipschitz function $\varphi$ defined on
$\Theta$ having finite $QH^1(\Theta)$ norm.
\item[2)] $W(uv)$ is defined as an element of
$L^\frac{2\omega}{\omega+1}(\Theta)$, with
$$\|W(uv)\|_{L^\frac{2\omega}{\omega+1}(\Theta)}\leq2C_4\|u\|_{QH^1(\Theta)}\|v\|_{QH^1(\Theta)},$$
with $C_4>0$ as in \eqref{globsob}, and with $W(uv)=uWv+vWu$ as elements
of $L^\frac{2\omega}{\omega+1}(\Theta)$.
\end{itemize}
\end{lem}

\textbf{Proof:} The proof of Lemma \ref{lem1} follows from the same
arguments used in \cite[Lemma 3.15]{ScottCJM}. Here, one
simply uses the global weak Poincar\'{e} inequality \eqref{globsob}
with gain $\omega>1$ in place of a global Sobolev inequality
as used in \cite{ScottCJM} (see condition \cite[(2.11)]{ScottCJM} or, equivalently, \eqref{GWP}). \begin{flushright}$\Box$\end{flushright}

\begin{defn}\label{def1}
Given a second order linear degenerate elliptic operator with rough coefficients $X$ as in \eqref{diffop}, we
introduce the associated bilinear form acting on $QH^1(\Theta)\times QH^1(\Theta)$,
\begin{equation}\label{bilform}
\mathcal{L}(u,v)=\int_{\Theta}\Big[\crochet{\nabla v}{P(x)\nabla
u}+v\mathbf{HR}u+u\mathbf{GS}v+Fuv\Big]\,dx.
\end{equation}
\end{defn}
As was done in \cite{ScottCJM} for the $X$-Dirichlet problem, the bilinear form (\ref{bilform}) will be used in a moment to define a notion of weak solution for the $X$-Neumann problem.  To begin the study of such objects, we show the boundedness of $\mathcal{L}$ on $QH^1(\Theta)$ followed by an almost-coercive estimate.

\begin{pro}\label{prop1}
Let $\Theta$ be a bounded domain with $\overline\Theta\subset\Omega$ and assume that the global weak
Poincar\'e inequality \eqref{globsob} with gain $\omega>1$
holds. Assume that $F\in L^t(\Theta)$ with
$t\geq\omega^\prime=\frac{\omega}{\omega-1}$ and  that
$\mathbf{G},\,\mathbf{H}\in \big[L^q(\Theta)\big]^N$ with
$q\geq2\omega^\prime$. Then there exists a constant
$C_6=C_6\big(C_1,C_4,N,\|F\|_{L^{\omega^\prime}(\Theta)},$\\
$\|\,|\mathbf{G}|\,\|_{L^{2\omega^\prime}(\Theta)}+|\,|\mathbf{H}|\,\|_{L^{2\omega^\prime}(\Theta)}\big)>0$
such that
\begin{equation}\label{5}
|\mathcal{L}(u,v)|\leq C_6\|u\|_{QH^1(\Theta)}\|v\|_{QH^1(\Theta)}
\end{equation}
for every $u,\,v\in QH^1(\Theta)$.
\end{pro}

\begin{cor}\label{cor1}
Under the hypotheses of Proposition \ref{prop1}, the bilinear form
$\mathcal{L}(\cdot,\cdot)$ introduced in Definition \ref{def1} is
well defined and continuous on $QH^1(\Theta)\times QH^1(\Theta)$.
\end{cor}

\textbf{Proof of Proposition \ref{prop1}:} Let $u,v\in QH^1(\Theta)$.  Then, by \eqref{4} and Schwarz's
inequality
\begin{eqnarray}
  \nonumber \left|\int_\Theta\crochet{\nabla v}{P(x)\nabla u}\,dx\right|&\leq& \int_\Theta|\sqrt{P(x)}\nabla
       u|\,|\sqrt{P(x)}\nabla v|\,dx\\
  \nonumber&\leq&\left(\int_\Theta|\sqrt{P(x)}\nabla u|^2\,dx\right)^\frac{1}{2}\left(\int_\Theta|\sqrt{P(x)}\nabla
       v|^2\,dx\right)^\frac{1}{2}\\
  \nonumber&\leq&C_1^2\|u\|_{QH^1(\Theta)}\|v\|_{QH^1(\Theta)}.
\end{eqnarray}
Using a generalization of H\"{o}lder's inequality (see
\cite[(7.11)]{GT} with exponents $2\omega$, $2\omega^\prime$, and
$2$) together with the results of Lemma \ref{lem1} we obtain
\begin{eqnarray}
  \nonumber \left|\int_\Theta v\mathbf{HR}u\,dx\right|&\leq&\|v\|_{L^{2\omega}(\Theta)}\|\,|\mathbf{H}|\,\|_{L^{2\omega^\prime}(\Theta)}
     \|\,|\mathbf{R}u|\,\|_{L^2(\Theta)}\\
  \nonumber&\leq&C_4\sqrt{N}\|\,|\mathbf{H}|\,\|_{L^{2\omega^\prime}(\Theta)}\|v\|_{QH^1(\Theta)}\|u\|_{QH^1(\Theta)}.
\end{eqnarray}
Similarly,
\begin{eqnarray}
  \nonumber \left|\int_\Theta u\mathbf{GS}v\,dx\right|&\leq&C_4\sqrt{N}\|\,|\mathbf{G}|\,\|_{L^{2\omega^\prime}(\Theta)}
     \|u\|_{QH^1(\Theta)}\|v\|_{QH^1(\Theta)}.
\end{eqnarray}
Another application of \cite[(7.11)]{GT} with exponents $\omega^\prime$, $2\omega$, and $2\omega$ gives us
\begin{eqnarray}
  \nonumber \left|\int_\Theta Fuv\,dx\right|&\leq&\|F\|_{L^{\omega^\prime}(\Theta)}\|u\|_{L^{2\omega}(\Theta)}
     \|v\|_{L^{2\omega}(\Theta)}\\
  \nonumber&\leq&C_4^2\|F\|_{L^{\omega^\prime}(\Theta)}\|u\|_{QH^1(\Theta)}\|v\|_{QH^1(\Theta)}.
\end{eqnarray}
Thus, \eqref{5} holds with
$$C_6=C_1^2+C_4\sqrt{N}\big(\|\,|\mathbf{G}|\,\|_{L^{2\omega^\prime}(\Theta)}+
\|\,|\mathbf{H}|\,\|_{L^{2\omega^\prime}(\Theta)}\big)+C_4^2\|F\|_{L^{\omega^\prime}(\Theta)}.$$
\begin{flushright}$\Box$\end{flushright}

\begin{pro}\label{prop2}
Let $\Theta$ be a bounded domain such that
$\overline\Theta\subset\Omega$ and assume that the global weak
Poincar\'{e} inequality \eqref{globsob} with gain $\omega>1$
holds. Assume that $F\in L^t(\Theta)$ with
$t>\omega^\prime$ and  that
$\mathbf{G},\,\mathbf{H}\in \big[L^q(\Theta)\big]^N$ with
$q>2\omega^\prime$. Then, there exists a constant
$C_7=C_7\big(c_1,C_4,\omega,N,t,q,\|F\|_{L^t(\Theta)},$\\
$\|\,|\mathbf{G}|\,\|_{L^q(\Theta)}+\|\,|\mathbf{H}|\,\|_{L^q(\Theta)}\big)>0$
so that
\begin{equation}\label{17}
|\mathcal{L}(u,u)|\geq
\frac{c_1}{4}\|u\|^2_{QH^1(\Theta)}-C_7\|u\|^2_{L^2(\Theta)}
\end{equation}
for every $u\in QH^1(\Theta)$.
\end{pro}

\textbf{Proof of Proposition \ref{prop2}:} Let $u\in
QH^1(\Theta)$.  Then, \eqref{4} gives
\begin{equation}\label{7}
\mathcal{L}(u,u)\geq c_1\int_{\Theta} \crochet{\nabla u}{Q(x)\nabla
u}\,dx-\int_\Theta|u(\mathbf{HR}u+\mathbf{GS}u)|\,dx-\int_\Theta|F|u^2\,dx.
\end{equation}
Using H\"{o}lder's inequality with exponents $t,t^\prime\geq1$ we
have
$$\int_\Theta|F|u^2\,dx\leq\|F\|_{L^t(\Theta)}\|u\|^2_{L^{2t^\prime}(\Theta)}\leq
\big(\e\|u\|^2_{L^{2\omega}(\Theta)}+\e^{-\frac{2\omega}{t(\omega-1)-\omega}}\|u\|^2_{L^2(\Theta)}\big)\|F\|_{L^t(\Theta)},$$
for any $\e>0$; we used the interpolation inequality \cite[(7.10)]{GT}. Now, if
$\|F\|_{L^t(\Theta)}=0$ then $\int_\Theta|F|u^2\,dx=0$, otherwise we
choose $\e=\frac{c_1}{8\|F\|_{L^t(\Theta)}C_4^2}$ so that \eqref{globsob} gives
\begin{equation}\label{8}
\int_\Theta|F|u^2\,dx\leq\alpha\|u\|^2_{L^2(\Theta)}+\frac{c_1}{8}\|u\|_{QH^1(\Theta)}^2,
\end{equation}
with $\alpha=\alpha(c_1,C_4,t,\omega,\|F\|_{L^t(\Theta)})>0$.

In a similar way, we use Lemma \ref{lem1} and \cite[(7.11)]{GT} to obtain
\begin{equation*}
\begin{array}{l}
 \displaystyle\int_\Theta|u(\mathbf{HR}u+\mathbf{GS}u)|\,dx\\
 \displaystyle\hspace{1cm}\leq\|u\|_{L^\frac{2q}{q-2}(\Theta)}
   \big(\|\,|\mathbf{H}|\,\|_{L^q(\Theta)}\|\,|\mathbf{R}u|\,\|_{L^2(\Theta)}+\|\,|\mathbf{G}|\,\|_{L^q(\Theta)}\|\,|\mathbf{S}u|\,\|_{L^2(\Theta)}\big)\\
 \displaystyle\hspace{1cm}\leq\sqrt{N}\big(\|\,|\mathbf{H}|\,\|_{L^q(\Theta)}+\|\,|\mathbf{G}|\,\|_{L^q(\Theta)}\big)
   \|u\|_{QH^1(\Theta)}\|u\|_{L^\frac{2q}{q-2}(\Theta)}\\
 \displaystyle\hspace{1cm}\leq\sqrt{N}\big(\|\,|\mathbf{H}|\,\|_{L^q(\Theta)}+\|\,|\mathbf{G}|\,\|_{L^q(\Theta)}\big)
   \|u\|_{QH^1(\Theta)}\big(\e\|u\|_{L^{2\omega}(\Theta)}+\e^{-\frac{2\omega}{(q-2)\omega-q}}\|u\|_{L^2(\Theta)}\big),
\end{array}
\end{equation*}
for any $\e>0$; we again used \cite[(7.10)]{GT}. We now use Young's inequality to show that for any $\delta>0$,
\begin{equation*}
\begin{array}{l}
 \displaystyle\int_\Theta|u(\mathbf{HR}u+\mathbf{GS}u)|\,dx\\
 \displaystyle\hspace{1cm}\leq\e\sqrt{N}\big(\|\,|\mathbf{H}|\,\|_{L^q(\Theta)}+\|\,|\mathbf{G}|\,\|_{L^q(\Theta)}\big)
   \|u\|_{QH^1(\Theta)}\|u\|_{L^{2\omega}(\Theta)}\\
 \displaystyle\hspace{2cm}+\e^{-\frac{2\omega}{(q-2)\omega-q}}\sqrt{N}\big(\|\,|\mathbf{H}|\,\|_{L^q(\Theta)}+\|\,|\mathbf{G}|\,\|_{L^q(\Theta)}\big)
   \|u\|_{QH^1(\Theta)}\|u\|_{L^2(\Theta)}\\
\displaystyle\hspace{1cm}\leq\e\sqrt{N}\big(\|\,|\mathbf{H}|\,\|_{L^q(\Theta)}+\|\,|\mathbf{G}|\,\|_{L^q(\Theta)}\big)
   \|u\|_{QH^1(\Theta)}\|u\|_{L^{2\omega}(\Theta)}\\
 \displaystyle\hspace{2cm}+\e^{-\frac{2\omega}{(q-2)\omega-q}}\sqrt{N}\big(\|\,|\mathbf{H}|\,\|_{L^q(\Theta)}+\|\,|\mathbf{G}|\,\|_{L^q(\Theta)}\big)
   \left(\frac{\delta}{2}\|u\|^2_{QH^1(\Theta)}+\frac{1}{2\delta}\|u\|_{L^2(\Theta)}^2\right).
\end{array}
\end{equation*}
If
$\|\,|\mathbf{H}|\,\|_{L^q(\Theta)}+\|\,|\mathbf{G}|\,\|_{L^q(\Theta)}=0$
then $\int_\Theta|u(\mathbf{HR}u+\mathbf{GS}u)|\,dx=0$, otherwise we
choose
$$\e=\frac{c_1}{8\sqrt{N}C_4\big(\|\,|\mathbf{H}|\,\|_{L^q(\Theta)}+\|\,|\mathbf{G}|\,\|_{L^q(\Theta)}\big)},\qquad
\delta=\frac{\e^{\frac{2\omega}{(q-2)\omega-q}}c_1}{\sqrt{N}\big(\|\,|\mathbf{H}|\,\|_{L^q(\Theta)}+\|\,|\mathbf{G}|\,\|_{L^q(\Theta)}\big)}$$
so that \eqref{globsob} gives
\begin{equation}\label{9}
\int_\Theta|u(\mathbf{HR}u+\mathbf{GS}u)|\,dx\leq\frac{5c_1}{8}\|u\|_{QH^1(\Theta)}^2+\beta\|u\|^2_{L^2(\Theta)},
\end{equation}
where
$\beta=\beta(c_1,C_4,N,q,\omega,\|\,|\mathbf{H}|\,\|_{L^q(\Theta)}+\|\,|\mathbf{G}|\,\|_{L^q(\Theta)})>0$.
Inserting \eqref{8} and \eqref{9} into \eqref{7} we have
$$\mathcal{L}(u,u)\geq\frac{c_1}{4}\|u\|^2_{QH^1(\Theta)}-\big(c_1+\alpha+\beta\big)\|u\|^2_{L^2(\Theta)},$$
which is \eqref{7} with $C_7=c_1+\alpha+\beta>0$.
\begin{flushright}$\Box$\end{flushright}

With the boundedness and almost-coercivity of the bilinear form $\mathcal{L}$ concluded, we now formally define the notion of weak solution associated to the $X$-Neumann problem.

\begin{defn}\label{def2}
Let $X$ be a second order linear degenerate elliptic operator with
rough coefficients as in \eqref{diffop}. Assume that
$\partial\Theta$ is piecewise $C^1$ and let $\nu$ be the unit
outward normal vector at each sufficiently regular boundary point.
Assume that the global weak Poincar\'e inequality (\ref{globsob})
with gain $\omega>1$ holds, let
$\mathbf{G},\,\mathbf{H}\in\big[L^q(\Theta)\big]^N$ with
$q\geq2\omega^\prime$ and let $F\in L^t(\Theta)$ with
$t\geq\omega^\prime$. If $f\in L^2(\Theta)$, $K\in\N$, $\mathbf{T}$
is a $K$--tuple of subunit vector fields and
$\mathbf{g}\in\big[L^2(\Theta)\big]^K$, then a function $(u,\nabla
u)\in QH^1(\Theta)$ is a weak solution of the $X$-Neumann Problem
\begin{equation}\label{NP}
\begin{cases}
 Xu=f+\mathbf{T^\prime}\mathbf{g}\qquad\qquad\qquad\qquad\quad\text{in }\Theta\\
 \crochet{\nu}{P(x)\nabla u+u\mathbf{G}\mathbf{S}-\mathbf{gT}}=0\,\qquad\text{on }\partial\Theta
\end{cases}
\end{equation}
if and only if
\begin{equation}\label{6}
\mathcal{L}(u,v)=\int_\Theta
fv+\mathbf{g}\mathbf{T}v\,dx\text{   $\quad$for every } v\in QH^1(\Theta).
\end{equation}
\end{defn}

\textbf{Motivation.} Assume that $\partial\Theta$ is smooth and that
the coefficients of the operator $X$ are $C^1(\overline\Theta)$. If
$u\in C^2(\overline\Theta)$ and $v\in C^1(\overline\Theta)$, it's
easy to see using the Divergence Theorem that $$\int_\Theta
v\,Xu\,dx=\mathcal{L}(u,v)-\int_{\partial\Theta}v\crochet{\nu}{P(x)\nabla
u+u\mathbf{GS}}\,d\sigma,$$ while
$$\int_{\Theta}\big(f+\mathbf{T^\prime g}\big)v=\int_{\Theta}fv+\mathbf{gT}v-\int_{\partial\Theta}v\crochet{\nu}{\mathbf{gT}}\,d\sigma.$$
Thus $u$ is a classical solution of problem \eqref{NP} if and only
if it satisfies \eqref{6}. Notice that if $\mathbf{G}$ and
$\mathbf{g}$ are identically null and if $P(x)$ is strictly positive
definite and continuous on $\overline\Theta$, we recover the
usual definition of weak solution of the Neumann Problem related to
the (now elliptic) differential operator $X$.\\

We now come to the first of our existence results - the existence of weak solutions to a modified version of the $X$-Neumann problem in $\Theta$.  The theorem is a consequence of Propositions (\ref{prop1}) and (\ref{prop2}) with an application of the Lax-Milgram lemma.

\begin{thm}\label{thm1}
Let $\Theta$ be a bounded domain such that
$\overline\Theta\subset\Omega$ and assume that the global weak
Poincar\'{e} inequality \eqref{globsob} with gain $\omega>1$ holds.
Let $X$ be a second order linear  degenerate elliptic operator with
rough coefficients as in \eqref{diffop}. Assume that $F\in
L^t(\Theta)$ with $t>\omega^\prime$ and that
$\mathbf{G},\,\mathbf{H}\in \big[L^q(\Theta)\big]^N$ with
$q>2\omega^\prime$. Then there is a constant $\gamma>0$ so that
given any $\mu\geq\gamma$, $f\in L^2(\Theta)$, $K\in\N$, any
$K$--tuple $\mathbf{T}$ of subunit vector fields and
$\mathbf{g}\in\big[L^2(\Theta)\big]^K$ there exists a unique weak
solution $u\in QH^1(\Theta)$ of the $X$-Neumann Problem
\begin{equation}\label{NPmu}
\begin{cases}
 Xu+\mu u=f+\mathbf{T^\prime}\mathbf{g}\qquad\qquad\,\,\,\,\quad\quad\text{in }\Theta\\
 \crochet{\nu}{P(x)\nabla u+u\mathbf{G}\mathbf{S}-\mathbf{gT}}=0\,\qquad\text{on
 }\partial\Theta.
\end{cases}
\end{equation}
Moreover, if $u\in QH^1(\Theta)$ is a weak solution of \eqref{NPmu}
then
\begin{equation}\label{solmap}
\|u\|_{QH^1(\Theta)}\leq\frac{4}{c_1}\big(\|f\|_{L^2(\Theta)}+\sqrt{K}\|\,|\mathbf{g}|\,\|_{L^2(\Theta)}\big).
\end{equation}
\end{thm}

\textbf{Proof of Theorem \ref{thm1}:} By Proposition \ref{prop2}
there is a $\gamma>0$ such that for every $\mu\geq\gamma$ one has
\begin{equation}\label{10}
\mathcal{L}(u,u)+\mu\|u\|^2_{L^2(\Theta)}\geq\mathcal{L}(u,u)+\gamma\|u\|^2_{L^2(\Theta)}\geq\frac{c_1}{4}\|u\|^2_{QH^1(\Theta)}
\end{equation}
for every $u\in QH^1(\Theta)$. Using Proposition \ref{prop1}, one
also sees that for each $\mu\geq\gamma$
$$\left|\mathcal{L}(u,v)+\mu\int_\Theta uv\,dx\right|\leq
(C_6+\mu)\|u\|_{QH^1(\Theta)}\|v\|_{QH^1(\Theta)}$$ for every
$u,v\in QH^1(\Theta)$. It follows that, for every
$\mu\geq\gamma$, the bilinear form $\mathcal{L}_\mu$ defined on
$QH^1(\Theta)\times QH^1(\Theta)$ by setting
$$\mathcal{L}_\mu(u,v)=\mathcal{L}(u,v)+\mu\int_\Theta uv\,dx\quad\text{for
every }u,v\in QH^1(\Theta)$$ is both bounded and coercive. Next, we
notice that Lemma \ref{lem1} implies that the map $\phi$ defined by $$\phi(v)=\int_\Theta
fv+\mathbf{g}\mathbf{T}v\,dx\quad\text{for every }v\in
QH^1(\Theta)$$ is linear and continuous. Indeed,
\begin{equation}\label{11}
\begin{array}{rl}
\displaystyle\left|\int_\Theta
fv+\mathbf{g}\mathbf{T}v\,dx\right|&\leq\|f\|_{L^2(\Theta)}\|v\|_{L^2(\Theta)}+\|\,|\mathbf{g}|\,\|_{L^2(\Theta)}\|\,|\mathbf{T}v|\,\|_{L^2(\Theta)}\\
\displaystyle&\leq\big(\|f\|_{L^2(\Theta)}+\sqrt{K}\|\,|\mathbf{g}|\,\|_{L^2(\Theta)}\big)\|v\|_{QH^1(\Theta)}
\end{array}
\end{equation}
for every $v\in QH^1(\Theta)$. Applying the Lax--Milgram lemma we conclude that there
exists a unique $u\in QH^1(\Theta)$ satisfying
$$\mathcal{L}_\mu(u,v)=\phi(v)$$
for every $v\in QH^1(\Theta)$.  Recalling Definition \ref{def2}, we see that $u$ is the unique weak solution in $QH^1(\Theta)$ of the $X$-Neumann problem \eqref{NPmu}.  Moreover, \eqref{10} and \eqref{11} indicate that such a $u$ also
satisfies
\bea\label{globest1}\frac{c_1}{4}\|u\|_{QH^1(\Theta)}\leq\big(\|f\|_{L^2(\Theta)}+\sqrt{K}\|\,|\mathbf{g}|\,\|_{L^2(\Theta)}\big)
\eea
and inequality \eqref{solmap} follows.
\begin{flushright}$\Box$\end{flushright}

\begin{rem}\label{rem2}
Notice that the same arguments used in the proof of Theorem
\ref{thm1} show that for every $\mu\geq\gamma$ and every
$\phi\in\big(QH^1(\Theta)\big)^*\cong QH^1(\Theta)$, the dual space to $QH^1(\Theta)$,
there exists a unique $u\in QH^1(\Theta)$ such that
\begin{equation}\label{12}
\mathcal{L}_\mu(u,v)=\mathcal{L}(u,v)+\mu\int_\Theta
uv\,dx=\phi(v)\quad\text{for every }v\in QH^1(\Theta).
\end{equation}
Moreover, the solution map $S_\mu:\big(QH^1(\Theta)\big)^*\rightarrow
QH^1(\Theta)$, defined by setting $S_\mu(\phi)=u$ if and only if
$u\in QH^1(\Theta)$ satisfies \eqref{12}, is linear and continuous,
with
$$\|u\|_{QH^1(\Theta)}\leq\frac{4}{c_1}\|\phi\|_{(QH^1(\Theta))^*}.$$
\end{rem}

\begin{defn}\label{def3}
Let $X$ be a second order linear degenerate elliptic operator with
rough coefficients as in \eqref{diffop}. Assume that
$\partial\Theta$ is piecewise $C^1$ and let $\nu$ be the outward
unit normal to $\partial\Theta$ at each sufficiently regular
boundary point. Assume that the global weak Poincar\'{e} inequality
(\ref{globsob}) with gain $\omega>1$ holds, let
$\mathbf{G},\,\mathbf{H}\in\big[L^q(\Theta)\big]^N$ with
$q\geq2\omega^\prime$ and let $F\in L^t(\Theta)$ with
$t\geq\omega^\prime$. If $f\in L^2(\Theta)$, $K\in\N$, $\mathbf{T}$
is a $K$--tuple of subunit vector fields and
$\mathbf{g}\in\big[L^2(\Theta)\big]^K$, then the adjoint problem to
the $X$-Neumann Problem \eqref{NP} is given by
\begin{equation}\label{NPadj}
\begin{cases}
 X^*u:=-\text{div}(P(x)\nabla
u)+\mathbf{G}\mathbf{S}u+\mathbf{R}^\prime(\mathbf{H}u)+Fu=f+\mathbf{T^\prime}\mathbf{g}\qquad\quad\text{in }\Theta\\
 \crochet{\nu}{P(x)\nabla u+u\mathbf{H}\mathbf{R}-\mathbf{gT}}=0\,\qquad\,\,\,\,\qquad\qquad\qquad\qquad\qquad\qquad\text{on
 }\partial\Theta.
\end{cases}
\end{equation}
We will say that the $X$-Neumann Problem \eqref{NP} is \emph{self--adjoint} if
$\mathbf{HR}\equiv\mathbf{GS}$ on $\Theta$, i.e. if $$\int_\Theta
v\mathbf{HR}u\,dx=\int_\Theta v\mathbf{GS}u\,dx\quad\text{for every
}u,v\in QH^1(\Theta).$$
\end{defn}

\begin{rem}\label{rem3}
The bilinear form $\mathcal{L}^*:QH^1(\Theta)\times
QH^1(\Theta)\rightarrow\R$ associated to the adjoint
problem \eqref{NPadj} is
\begin{equation}\label{14}
\mathcal{L}^*(u,v)=\int_{\Theta}\crochet{\nabla v}{P(x)\nabla
u}+v\mathbf{GS}u+u\mathbf{HR}v+Fuv\,dx=\mathcal{L}(v,u),
\end{equation}
for every $u,v\in QH^1(\Theta)$.
\end{rem}

We now present our main existence result for the $X$-Neumann problem.

\begin{thm}\label{thm2}
Let $(\Omega,\rho)$ be a geometric homogeneous space and $\Theta$ be a bounded domain with $\overline\Theta\subset\Omega$. Assume that that the local Poincar\'e inequality (\ref{Poincare-loc}) holds with $p=2$ and that the global weak
Poincar\'{e} inequality \eqref{globsob} with gain $\omega>1$
holds. Let $X$ be a
second order linear degenerate elliptic operator with
rough coefficients as in \eqref{diffop}. Assume that $F\in
L^t(\Theta)$ with $t>\omega^\prime$ and
that $\mathbf{G},\,\mathbf{H}\in \big[L^q(\Theta)\big]^N$ with
$q>2\omega^\prime$. Then, we have the following conclusions.
\begin{itemize}
\item[1)] One and only one of the following alternatives hold:\\
{\bf either}
\begin{itemize}
\item[(I)] given any $f\in L^2(\Theta)$, $K\in\N$, any $K$--tuple of subunit vector fields $\mathbf{T}$, and
$\mathbf{g}\in\big[L^2(\Theta)\big]^K$ there exists a unique weak
solution $u\in QH^1(\Theta)$ of the $X$-Neumann Problem \eqref{NP}
\end{itemize}
{\bf or}
\begin{itemize}

\item[(II)] there exist nontrivial weak solutions $u\in
QH^1(\Theta)\setminus\{(0,{\bf 0})\}$ of the homogeneous $X$-Neumann
Problem
\begin{equation}\label{NPhom}
\begin{cases}
 Xu=0\qquad\qquad\qquad\qquad\,\,\quad\text{in }\Theta\\
 \crochet{\nu}{P(x)\nabla u+u\mathbf{G}\mathbf{S}}=0\,\qquad\text{on
 }\partial\Theta.
\end{cases}
\end{equation}
\end{itemize}
\item[2)] If alternative 1)--(II) holds, then the dimension of the subspace
$N\subset QH^1(\Theta)$ of weak solutions of the homogeneous $X$-Neumann
Problem \eqref{NPhom} is finite and equals the dimension of the
subspace $N^*\subset QH^1(\Theta)$ of weak solutions of the adjoint
homogeneous $X$-Neumann Problem
\begin{equation}\label{NPadjhom}
\begin{cases}
 X^*w=0\qquad\qquad\qquad\qquad\quad\text{in }\Theta\\
 \crochet{\nu}{P(x)\nabla w+w\mathbf{HR}}=0\,\qquad\text{on
 }\partial\Theta.
\end{cases}
\end{equation}
\item[3)] If alternative 1)--(II) holds, the $X$-Neumann Problem \eqref{NP}
admits a weak solution $u\in QH^1(\Theta)$ for given $f\in
L^2(\Theta)$, $K\in\N$, $\mathbf{T}$ a $K$--tuple of subunit vector
fields and $\mathbf{g}\in\big[L^2(\Theta)\big]^K$ if and only if
$$\int_\Theta fw+\mathbf{gT}w\,dx=0$$
for every $w\in N^*$.
\end{itemize}
\end{thm}

\textbf{Proof of Theorem \ref{thm2}:} To begin, we notice that Theorem \ref{thm1} and Remark \ref{rem3} provide a $\gamma>0$
so that for every $\mu\geq\gamma$ and every $\phi\in\big(QH^1(\Theta)\big)^*$, there exist unique $v,w\in QH^1(\Theta)$ so that
$$\mathcal{L}(v,\psi)+\mu\int_\Theta v\psi\,dx=\phi(\psi),\quad\text{and}\quad\mathcal{L}^*(w,\psi)+\mu\int_\Theta w\psi\,dx=\phi(\psi)$$
for every $\psi\in QH^1(\Theta)$. The corresponding solution maps
$$S_\mu,S^*_\mu:\big(QH^1(\Theta)\big)^*\rightarrow QH^1(\Theta)$$
defined by setting $S_\mu\phi=v$ and $S^*_\mu\phi=w$ are well
defined, linear and continuous. Next, we introduce the map
$J:QH^1(\Theta)\rightarrow\big(QH^1(\Theta)\big)^*$ defined by
setting $$Ju(v)=\int_\Theta uv\,dx$$ for every $u,v\in
QH^1(\Theta)$. The map $J$ is linear, continuous and, by Proposition
\ref{compact2} with $q=2$, compact.

Now, given $f\in L^2(\Theta)$, $K\in\N$, $\mathbf{T}$ a $K$--tuple
of subunit vector fields and $\mathbf{g}\in\big[L^2(\Theta)\big]^K$
we define a continuous linear functional $\vphi$ on $QH^1(\Theta)$:
\begin{equation}\label{16}
\vphi(\psi)=\int_\Theta f\psi+\mathbf{gT}\psi\,dx.
\end{equation}
Notice that $u\in QH^1(\Theta)$ is a weak solution of the $X$-Neumann
Problem \eqref{NP} if and only if
\begin{equation*}
\begin{cases}
 Xu+\gamma u=\gamma u+f+\mathbf{T^\prime}\mathbf{g}\qquad\,\,\quad\text{in }\Theta\\
 \crochet{\nu}{P(x)\nabla u+u\mathbf{G}\mathbf{S}-\mathbf{gT}}=0\,\quad\text{on }\partial\Theta;
\end{cases}
\end{equation*}
that is, if and only if one has
$$\mathcal{L}(u,\psi)+\gamma\int_\Theta u\psi\,dx=\vphi(\psi)+\gamma Ju(\psi)$$
for every $\psi\in QH^1(\Theta)$.  This is equivalent to requiring that $u\in QH^1(\Theta)$ solves
\begin{equation}\label{13}
u-\mathcal{K}u=\Phi
\end{equation}
with $\Phi=S_\gamma\vphi\in QH^1(\Theta)$ where $\mathcal{K}=\gamma S_\gamma\circ J:QH^1(\Theta)\rightarrow
QH^1(\Theta)$; a compact linear operator. Similarly we see that $u\in
QH^1(\Theta)$ solves the homogeneous $X$-Neumann Problem \eqref{NPhom}
if and only if it is a solution of $$u-\mathcal{K}u=0.$$ Hence, we
can apply the Fredholm Alternative, see for instance \cite[theorem 5 (appendix D)]{Evans}, and conclude that either
\begin{itemize}
\item[(A)] for every $\Phi\in QH^1(\Theta)$ equation \eqref{13}
admits a unique solution $u\in QH^1(\Theta)$,
\end{itemize}
or
\begin{itemize}
\item[(B)] the homogeneous equation associated to \eqref{13}
admits nontrivial solutions $u\in QH^1(\Theta)\setminus\{(0,{\bf 0})\}$.
\end{itemize}

If alternative (A) holds, then clearly the $X$-Neumann Problem
\eqref{NP} admits a unique weak solution in $QH^1(\Theta)$ for any
choice of $f\in L^2(\Theta)$, $K\in\N$, $\mathbf{T}$ a $K$--tuple of
subunit vector fields and $\mathbf{g}\in\big[L^2(\Theta)\big]^K$;
this proves 1)--(I). On the other hand, if alternative (B) holds,
the Fredholm Alternative states that there are nontrivial solutions
$u\in QH^1(\Theta)\setminus\{(0,{\bf 0})\}$ of the homogeneous
equation $u-\mathcal{K}u=0$. This proves 1)--(II) and completes the
proof of 1).

Assume now that alternative (B) holds, i.e. that 1)--(II) of the statement of the Theorem holds.  Let
$\mathcal{K}^*:\big(QH^1(\Theta)\big)^*\rightarrow\big(QH^1(\Theta)\big)^*$
be the adjoint operator to $\mathcal{K}$ and let $V^*\subset\big(QH^1(\Theta)\big)^*$ be the subspace of weak solutions of
\begin{equation}\label{15}
\Psi-\mathcal{K}^*\Psi=0.
\end{equation}
Then, by the Fredholm Alternative we obtain
$$1\leq\text{dim}N=\text{dim}V^*<\infty.$$
By the properties of adjoint operators one has
$$\mathcal{K}^*=\gamma J^*\circ(S_\gamma)^*,$$ where $J^*$,
$(S_\gamma)^*$ are the adjoint operators of $J$ and
$S_\gamma$ respectively. Since
$QH^1(\Theta)$ is reflexive and $\mathcal{L}(u,v)=\mathcal{L}^*(v,u)$ for every $u,v\in
QH^1(\Theta)$, it is evident that $J^*\equiv J$ and $(S_\gamma)^*\equiv
S^*_\gamma$. Therefore, we have that $\Psi\in\big(QH^1(\Theta)\big)^*$ is
a solution of equation \eqref{15} if and only if
$$\Psi=\mathcal{K}^*\Psi=\gamma J\big(S^*_\gamma\Psi\big).$$ Setting
$w=S^*_\gamma\Psi\in QH^1(\Theta)$ and recalling the definitions of
the mappings $J$ and $S^*_\gamma$, we obtain
$$\Psi(v)=\gamma\int_\Theta vw\,dx\quad\text{ for every }v\in
QH^1(\Theta),$$ where $w$ is a weak solution of the adjoint
homogeneous $X$-Neumann Problem \eqref{NPadjhom}. This implies that
$V^*=JN^*$, so that
$$\text{dim}V^*\leq\text{dim}N^*$$ as the map $J$ may fail to be
injective (we recall that the map $i: QH^1(\Theta)\ra L^2(\Theta)$
of Proposition (\ref{compact2}) may not be injective). We now
conclude that
$$\text{dim}N\leq\text{dim}N^*.$$  If we repeat the above
argument, replacing the operator $X$ with $X^*$ and the bilinear
form $\mathcal{L}(\cdot,\cdot)$ with $\mathcal{L}^*(\cdot,\cdot)$,
we arrive at the opposite inequality:
$\text{dim}N\geq\text{dim}N^*$. Part 2) of the theorem had now been proven.

Suppose again that alternative (B) holds. That is, assume that item 1)--(II) holds. By the previous
arguments, $u\in QH^1(\Theta)$ is a weak solution of the $X$-Neumann
Problem \eqref{NP} if and only if it is a solution of equation
\eqref{13}, with $\Phi=S_\gamma\vphi\in QH^1(\Theta)$ and where
$\vphi\in\big(QH^1(\Theta)\big)^*$ is defined by \eqref{16}. The Fredholm Alternative indicates that such an equation admits a solution if and only if
$\Psi(\Phi)=0$ for every $\Psi\in V^*=JN^*$. That is, if and only if
$$\int_\Theta\Phi w\,dx=0$$ for every weak solution $w\in QH^1(\Theta)$ of the homogeneous adjoint $X$-Neumann Problem
\eqref{NPadjhom}. Notice that such $w$ satisfy
$$\int_\Theta fw+\mathbf{gT}w\,dx\,=\,\varphi(w)\,=\,\mathcal{L}_\gamma(\Phi,w)\,=\,\mathcal{L}^*(w,\Phi)+\gamma\int_\Theta \Phi w\,dx\,
=\,\gamma\int_\Theta \Phi w\,dx.$$ Thus, problem
\eqref{NP} admits a weak solution $u\in QH^1(\Theta)$ if and only if
$$\int_\Theta fw+\mathbf{gT}w\,dx=0$$
for every $w\in N^*$.  This completes the proof of part 3).
\begin{flushright}$\Box$\end{flushright}

With Theorem \ref{thm2} in hand we now begin our analysis of
spectral properties associated to the $X$-Neumann problem.

\begin{thm}\label{thm3}
Let $(\Omega,\rho)$ be a geometric homogeneous space and let $\Theta$ be a bounded domain with $\overline\Theta\subset\Omega$. Assume that both the local Poincar\'e inequality (\ref{Poincare-loc}) with $p=2$ and the global weak
Poincar\'{e} inequality \eqref{globsob} with gain $\omega>1$
hold.  Let $X$ be a second order linear degenerate elliptic operator with
rough coefficients as in \eqref{diffop}. Assume that $F\in
L^t(\Theta)$ with $t>\omega^\prime$ and
that $\mathbf{G},\,\mathbf{H}\in \big[L^q(\Theta)\big]^N$ with
$q>2\omega^\prime$. Then there exists an at most countable set
$\Sigma\subset\R$ so that the problem
\begin{equation}\label{NPeigen}
\begin{cases}
 Xu=\lambda u+f+\mathbf{T^\prime}\mathbf{g}\qquad\qquad\qquad\quad\text{in }\Theta\\
 \crochet{\nu}{P(x)\nabla u+u\mathbf{G}\mathbf{S}-\mathbf{gT}}=0\,\,\,\qquad\text{on
 }\partial\Theta
\end{cases}
\end{equation}
has a unique weak solution $u\in QH^1(\Theta)$ for every $f\in
L^2(\Theta)$, every $K\in\N$, every $K$--tuple $\mathbf{T}$ of
subunit vector fields and every
$\mathbf{g}\in\big[L^2(\Theta)\big]^K$ if and only if
$\lambda\notin\Sigma$. Moreover, if $\Sigma$ is infinite, its
elements can be arranged in a monotone sequence that diverges to
$+\infty$.
\end{thm}

\begin{defn}
With $\Sigma$ as in Theorem \ref{thm3}, we will call any $\lambda\in\Sigma$ an eigenvalue of the $X$-Neumann
Problem \eqref{NPeigen}. Any weak solution $u\in
QH^1(\Theta)\setminus\{(0,{\bf 0})\}$ of the Homogeneous $X$-Neumann Problem
\begin{equation}\label{NPeigenhom}
\begin{cases}
 Xu=\lambda u\qquad\qquad\qquad\qquad\,\,\,\,\quad\quad\text{in }\Theta\\
 \crochet{\nu}{P(x)\nabla u+u\mathbf{G}\mathbf{S}}=0\qquad\qquad\text{on
 }\partial\Theta
\end{cases}
\end{equation}
will be called an eigenfunction of the $X$-Neumann Problem \eqref{NPeigen} associated to the eigenvalue $\lambda$.
\end{defn}

\textbf{Proof of Theorem \ref{thm3}:} By Theorem \ref{thm1} there
exists a $\gamma>0$ such that for every $\mu\geq\gamma$, every $f\in
L^2(\Theta)$, every $K\in\N$, every $K$--tuple $\mathbf{T}$ of
subunit vector fields and every
$\mathbf{g}\in\big[L^2(\Theta)\big]^K$ there exists a unique weak
solution to problem \eqref{NPmu}. Thus, problem \eqref{NPeigen}
admits a unique weak solution whenever $\lambda\leq-\gamma$.

From now on we will assume that $\lambda>-\gamma$, with $\gamma>0$.
Arguing as in the proof of Theorem \ref{thm2} we see that problem
\eqref{NPeigen} admits a unique weak solution for every $f\in
L^2(\Theta)$, every $K\in\N$, every $K$--tuple $\mathbf{T}$ of
subunit vector fields and every
$\mathbf{g}\in\big[L^2(\Theta)\big]^K$ if and only if $(u,\nabla
u)=(0,\mathbf{0})$ is the unique weak solution of
\eqref{NPeigenhom}. Moreover, if problem \eqref{NPeigenhom} admits
nontrivial weak solutions, then the subspace generated by those weak
solutions has finite dimension.

Now, $u\in QH^1(\Theta)$ is a weak solution of \eqref{NPeigenhom} if
and only if it is a weak solution of
\begin{equation*}
\begin{cases}
 Xu+\gamma u=(\lambda+\gamma) u\qquad\qquad\qquad\quad\text{in }\Theta\\
 \crochet{\nu}{P(x)\nabla u+u\mathbf{G}\mathbf{S}-\mathbf{gT}}=0\,\,\,\,\qquad\text{on
 }\partial\Theta,
\end{cases}
\end{equation*}
and this in turn holds if and only if $$u=S_\gamma\circ
J\big((\gamma+\lambda)u\big)=\frac{\gamma+\lambda}{\gamma}\mathcal{K}u,$$
where $\mathcal{K}=\gamma S_\gamma\circ J:QH^1(\Theta)\rightarrow
QH^1(\Theta)$ is the linear, compact operator defined in the proof of Theorem \ref{thm2}. Thus,
$u\in QH^1(\Theta)\setminus\{(0,{\bf 0})\}$ is a weak solution of
problem \eqref{NPeigenhom} if and only if it is an eigenfunction of
the compact linear operator $\mathcal{K}$ associated to the
eigenvalue $\frac{\gamma}{\lambda+\gamma}$.

The set $\Sigma^\prime$ of real eigenvalues of $\mathcal{K}$ is
countable at most and, if it is infinite, its elements can be
arranged as a sequence converging to $0$. Consequently, the set
$\Sigma\subset\R$ of numbers $\lambda$ such that problem
\eqref{NPeigenhom} has nontrivial weak solutions in
$QH^1(\Theta)\setminus\{(0,{\bf 0})\}$ is countable at most and, if infinite,
it comprises the values of a monotone sequence diverging to
$+\infty$.

 \begin{flushright}$\Box$\end{flushright}

\begin{thm}\label{thm4}
Let $(\Omega,\rho)$ be a geometric homogeneous space and let
$\Theta$ be a bounded domain with $\overline\Theta\subset\Omega$.
Assume that both the local Poincar\'e inequality
(\ref{Poincare-loc}) with $p=2$ and the global weak Poincar\'{e}
inequality \eqref{globsob} with gain $\omega>1$ hold.  Let $X$ be a
second order linear degenerate elliptic operator with rough
coefficients as in \eqref{diffop}. Assume that $F\in L^t(\Theta)$
with $t>\omega^\prime$ and that $\mathbf{G},\,\mathbf{H}\in
\big[L^q(\Theta)\big]^N$ with $q>2\omega^\prime$. If
$\lambda\notin\Sigma$, there exists a positive constant
$C_8=C_8(\lambda,\Theta,\Omega,c_1,C_1,\mathbf{G},\mathbf{H},F)$
such that if $f\in L^2(\Theta)$, $K\in\N$, $\mathbf{T}$ is a
$K$--tuple of subunit vector fields,
$\mathbf{g}\in\big[L^2(\Theta)\big]^K$ and $u\in QH^1(\Theta)$ is
the unique weak solution of problem \eqref{NPeigen}, then one has
the estimate
\begin{equation}\label{22}
\|u\|_{QH^1(\Theta)}\leq
C_8\big(\|f\|_{L^2(\Theta)}+\sqrt{K}\|\,|\mathbf{g}|\,\|_{L^2(\Theta)}\big).
\end{equation}
\end{thm}

\textbf{Proof of Theorem \ref{thm4}:} We start by showing that under
the current hypotheses there exists a constant $\hat{C}>0$, independent of
$u,\,f,\,K,\,\mathbf{T},\,\mathbf{g}$, such that
\begin{equation}\label{21}
\|u\|_{L^2(\Theta)}\leq
\hat{C}\big(\|f\|_{L^2(\Theta)}+\sqrt{K}\|\,|\mathbf{g}|\,\|_{L^2(\Theta)}\big).
\end{equation}
To arrive at a contradiction, suppose that (\ref{21}) is false.
Then, for every $n\in\N$ there exist $f_n\in L^2(\Theta)$,
$K_n\in\N$, $\mathbf{g}_n\in\big[L^2(\Theta)\big]^{K_n}$,
$\mathbf{T}_n$ a $K_n$--tuple of subunit vector fields and $u_n \in
QH^1(\Theta)$ such that
\begin{equation*}
\begin{cases}
 Xu_n=\lambda u_n+f_n+\mathbf{T}_n^\prime\mathbf{g}_n\qquad\qquad\qquad\quad\text{in }\Theta\\
 \crochet{\nu}{P(x)\nabla u_n+u_n\mathbf{G}\mathbf{S}-\mathbf{g}_n\mathbf{T}_n}=0\,\,\,\,\qquad\text{on
 }\partial\Theta,
\end{cases}
\end{equation*}
that is
\begin{equation}\label{19}
\mathcal{L}(u_n,v)=\int_\Theta\lambda
u_nv+f_nv+\mathbf{g}_n\mathbf{T}_nv\,dx
\end{equation}
for every $v\in QH^1(\Theta)$, and
$$\|u_n\|_{L^2(\Theta)}>
n\big(\|f_n\|_{L^2(\Theta)}+\sqrt{K_n}\|\,|\mathbf{g}_n|\,\|_{L^2(\Theta)}\big).$$
Without loss of generality we can assume that
$\|u_n\|_{L^2(\Theta)}=1$ for every $n\in\N$, so that
\begin{equation}\label{18}
\big(\|f_n\|_{L^2(\Theta)}+\sqrt{K_n}\|\,|\mathbf{g}_n|\,\|_{L^2(\Theta)}\big)<\frac{1}{n}
\end{equation}
for each $n\in\mathbb{N}$. Let $\gamma>0$ be as in Theorem \ref{thm1}. Since $u_n\in
QH^1(\Theta)$ is also a weak solution of
\begin{equation*}
\begin{cases}
 Xu_n+\gamma u_n=(\gamma+\lambda) u_n+f_n+\mathbf{T}_n^\prime\mathbf{g}_n\qquad\quad\quad\text{in }\Theta\\
 \crochet{\nu}{P(x)\nabla u_n+u_n\mathbf{G}\mathbf{S}-\mathbf{g}_n\mathbf{T}_n}=0\qquad\,\,\,\,\,\,\qquad\text{on
 }\partial\Theta,
\end{cases}
\end{equation*}
we obtain by inequality \eqref{solmap} that
\begin{eqnarray}
\nonumber
\|u_n\|_{QH^1(\Theta)}&\leq&\frac{4}{c_1}\big(\|(\gamma+\lambda)u_n+f_n\|_{L^2(\Theta)}+\sqrt{K_n}\|\,|\mathbf{g}_n|\,\|_{L^2(\Theta)}\big)\\
\nonumber&\leq&\frac{4}{c_1}\big((\gamma+\lambda)+1\big)
\end{eqnarray}
for every $n\in\N$. Since $QH^1(\Theta)$ is a Hilbert space and
since $QH^1(\Theta)$ is compactly embedded in $L^2(\Theta)$ by Proposition \ref{compact2}, we can assume (up to a subsequence) that
\begin{equation}\label{20}
\begin{array}{rcl}
(u_n,\nabla u_n)\!&\!\rightharpoonup\!&\!(u,\nabla u)\quad\text{in
}QH^1(\Theta),\text{ and}\\
u_n\!&\!\rightarrow\!&\!u\;\;\;\;\;\;\quad\quad\text{in }L^2(\Theta).
\end{array}
\end{equation}
Passing to the limit in \eqref{19}
while exploiting \eqref{18}, \eqref{20}, and the continuity of the
bilinear form $\mathcal{L}(\cdot,\cdot)$, we see that for every
$v\in QH^1(\Theta)$ we have
$$\mathcal{L}(u,v)=\int_\Theta\lambda uv\,dx,$$
since
\begin{eqnarray}
\nonumber\left|\int_\Theta
f_nv+\mathbf{g}_n\mathbf{T}_nv\,dx\right|\!&\!\leq\!&\!
\|f_n\|_{L^2(\Theta)}\|v\|_{L^2(\Theta)}+\|\,|\mathbf{g}_n|\,\|_{L^2(\Theta)}\|\,|\mathbf{T}_nv|\,\|_{L^2(\Theta)}\\
\nonumber\!&\!\leq\!&\!\big(\|f_n\|_{L^2(\Theta)}+\sqrt{K_n}\|\,|\mathbf{g}_n|\,\|_{L^2(\Theta)}\big)\|v\|_{QH^1(\Theta)}\,\,<\,\,\frac{1}{n}\|v\|_{QH^1(\Theta)}
\end{eqnarray}
for every $v\in QH^1(\Theta)$. Hence, $u\in QH^1(\Theta)$ is a weak
solution of problem (\ref{NPeigenhom})
and, since $\lambda\notin\Sigma$, we conclude that $(u,\nabla u)=(0,\mathbf{0})$. This is in contradiction with
$$\|u\|_{L^2(\Theta)}=\lim_{n\rightarrow\infty}\|u_n\|_{L^2(\Theta)}=1$$
and, therefore, establishes inequality \eqref{21}.\\

Now, if $u\in QH^1(\Theta)$ is a weak solution of \eqref{NPeigen} it
is also a weak solution of
\begin{equation*}
\begin{cases}
 Xu+\gamma u=(\lambda+\gamma) u+f+\mathbf{T^\prime}\mathbf{g}\qquad\quad\text{in }\Theta\\
 \crochet{\nu}{P(x)\nabla u+u\mathbf{G}\mathbf{S}-\mathbf{gT}}=0\,\,\,\quad\qquad\text{on
 }\partial\Theta,
\end{cases}
\end{equation*}
where $\gamma>0$ is as in Theorem \ref{thm1}. By inequalities
\eqref{solmap} and \eqref{21} we obtain
\begin{eqnarray}
\nonumber
\|u\|_{QH^1(\Theta)}&\leq&\frac{4}{c_1}\big(\|(\gamma+\lambda)u+f\|_{L^2(\Theta)}+\sqrt{K}\|\,|\mathbf{g}|\,\|_{L^2(\Theta)}\big)\\
\nonumber&\leq&\frac{4}{c_1}\big((\gamma+\lambda)\|u\|_{L^2(\Theta)}+\|f\|_{L^2(\Theta)}+\sqrt{K}\|\,|\mathbf{g}|\,\|_{L^2(\Theta)}\big)\\
\nonumber&\leq&\frac{4\big((\gamma+\lambda)\hat{C}+1\big)}{c_1}\big(\|f\|_{L^2(\Theta)}+\sqrt{K}\|\,|\mathbf{g}|\,\|_{L^2(\Theta)}\big),
\end{eqnarray}
which is inequality \eqref{22} with
$C_8:=\frac{4\left((\gamma+\lambda)\hat{C}+1\right)}{c_1}$.
\begin{flushright}$\Box$\end{flushright}

\begin{thm}\label{thm5}
Let $(\Omega,\rho)$ be a geometric homogeneous space and let $\Theta$ be a bounded domain such that
$\overline\Theta\subset\Omega$. Assume that both the local Poincar\'e inequality \eqref{Poincare-loc} with $p=2$ and the global weak
Poincar\'{e} inequality \eqref{globsob} with gain $\omega>1$
hold. Let $X$ be a
second order linear degenerate elliptic operator with
rough coefficients as in \eqref{diffop} that satisfies negativity condition (1) as in Definition \ref{neg1}. Assume that $F\in
L^t(\Theta)$ with $t>\omega^\prime$ and
that $\mathbf{G},\,\mathbf{H}\in \big[L^q(\Theta)\big]^N$ with
$q>2\omega^\prime$. Then the only weak solution  $(u,\nabla u)\in
QH^1(\Theta)$ of the $X$-Neumann Problem
\begin{equation}\label{23}
\begin{cases}
 Xu=0\qquad\qquad\qquad\qquad\quad\,\,\text{in }\Theta\\
 \crochet{\nu}{P(x)\nabla u+u\mathbf{G}\mathbf{S}}=0\,\qquad\text{on }\partial\Theta
\end{cases}
\end{equation}
is $(u,\nabla u)=(0,\mathbf{0})$.
\end{thm}

\textbf{Proof of Theorem \ref{thm5}:} We start by noticing that if
$\mathbf{G}=\mathbf{H}=0$ almost everywhere in $\Theta$,
$$0=\mathcal{L}(u,u)\geq c_1\int_\Theta\crochet{\nabla u}{Q(x)\nabla u}\,dx+\e\int_\Theta u^2\,dx\geq\min\{\e,c_1\}\|u\|^2_{QH^1(\Theta)}$$
where we have used that $X$ satisfies negativity condition (1).  Thus, in this case, we see that $(u,\nabla u)=(0,\mathbf{0})$. From this point onward, we shall assume that
$\|\,|\mathbf{G}|+|\mathbf{H}|\,\|_{L^q(\Theta)}\neq0$.

We now proceed by assuming that there is an $\e>0$ such that for
every $(u,\nabla u),(v,\nabla v)\in QH^1(\Theta)$ satisfying
$uv\geq0$ almost everywhere in $\Theta$ one has
\begin{equation}\label{24}
\int_\Theta Fuv+\mathbf{GS}(uv)\,dx\geq\e\int_\Theta uv\,dx,
\end{equation}
see Definition \ref{neg1}. Let $(u,\nabla u)\in
QH^1(\Theta)$ be a weak solution of \eqref{23}. Then for each $k>0$,
 $(v,\nabla v)=\big((u-k)^+,\chi_{\{u>k\}}\nabla
u\big)\in QH^1(\Theta)$ is a valid test function for $u$, see \cite{SW2}. Thus,
$$0=\mathcal{L}(u,v)=\int_\Theta\crochet{\nabla v}{P(x)\nabla u}\,dx+\int_\Theta v\mathbf{HR}u\,dx+\int_\Theta u\mathbf{GS}v\,dx+\int_\Theta Fuv\,dx.$$
Using Lemma \ref{lem1} and that $uv\geq0$
almost everywhere in $\Theta$, \eqref{24} gives us:
\begin{eqnarray}
\nonumber\int_\Theta\crochet{\nabla v}{P(x)\nabla
v}\,dx\!&\!=\!&\!-\int_\Theta\mathbf{GS}(uv)\,dx+\int_\Theta
v\mathbf{GS}u\,dx-\int_\Theta v\mathbf{HR}u\,dx-\int_\Theta
Fuv\,dx\\
\nonumber\!&\!\leq\!&\!-\e\int_\Theta uv\,dx+\int_\Theta
v(\mathbf{GS}u-\mathbf{HR}u)\,dx\\
\nonumber\!&\!\leq\!&\!-\e\int_\Theta v^2-\e k\int_\Theta
v\,dx+\int_\Theta v(\mathbf{GS}v-\mathbf{HR}v)\,dx.
\end{eqnarray}
Setting $\Gamma=\text{supp}\big(|\sqrt{Q(x)}\nabla v|\big)$
and $\e_0=\min\{c_1,\e\}$, we obtain
\begin{eqnarray}
\nonumber\e_0\|v\|^2_{QH^1(\Theta)}\!&\!=\!&\!\e_0\left(\int_\Theta|\sqrt{Q(x)}\nabla
v|^2\,dx+\int_\Theta v^2\,dx\right)\\
\nonumber\!&\!\leq\!&\!\int_\Theta|\sqrt{P(x)}\nabla
v|^2\,dx+\e\int_\Theta v^2\,dx\\
\nonumber\!&\!\leq\!&\!\int_\Theta|v|\big[|\mathbf{H}||\mathbf{R}v|+|\mathbf{G}||\mathbf{S}v|\big]\,dx.
\end{eqnarray}
Since \cite[lemma 3.18]{ScottCJM} shows that
$$|\mathbf{R}v|\leq\sqrt{N}|\sqrt{Q(x)}\nabla v|\qquad\text{and}\qquad|\mathbf{S}v|\leq\sqrt{N}|\sqrt{Q(x)}\nabla v|$$
almost everywhere in $\Theta$, we have
\begin{eqnarray}
\nonumber\e_0\|v\|^2_{QH^1(\Theta)}\!&\!\leq\!&\!\int_\Gamma|v|\big[|\mathbf{H}||\mathbf{R}v|+|\mathbf{G}||\mathbf{S}v|\big]\,dx\\
\nonumber\!&\!\leq\!&\!\sqrt{N}\big\|\big(|\mathbf{G}|+|\mathbf{H}|\big)v\big\|_{L^2(\Gamma)}\|\,|\sqrt{Q(x)}\nabla
v|\,\|_{L^2(\Gamma)}\\
\nonumber\!&\!\leq\!&\!\sqrt{N}\big\|\big(|\mathbf{G}|+|\mathbf{H}|\big)v\big\|_{L^2(\Gamma)}\|v\|_{QH^1(\Theta)},
\end{eqnarray}
and so H\"older's inequality and \eqref{Poincare-loc} give
\begin{eqnarray}
\nonumber\|v\|_{QH^1(\Theta)}\!&\!\leq\!&\!\frac{\sqrt{N}}{\e_0}\|\,|\mathbf{G}|+|\mathbf{H}|\,\|_{L^q(\Theta)}\|v\|_{L^\frac{2q}{q-2}(\Gamma)}\\
\label{25}\!&\!\leq\!&\!\frac{\sqrt{N}}{\e_0}\|\,|\mathbf{G}|+|\mathbf{H}|\,\|_{L^q(\Theta)}|\Gamma|^{\frac{1}{2\omega^\prime}-\frac{1}{q}}
  \|v\|_{L^{2\omega}(\Theta)}\\
\nonumber\!&\!\leq\!&\!C_3\frac{\sqrt{N}}{\e_0}\|\,|\mathbf{G}|+|\mathbf{H}|\,\|_{L^q(\Theta)}|\Gamma|^{\frac{1}{2\omega^\prime}-\frac{1}{q}}
  \|v\|_{QH^1(\Theta)}.
\end{eqnarray}
Dividing (\ref{25}) by $||v||_{QH^1(\Theta)}$ we obtain
\begin{equation}\label{26}
|\Gamma|\geq\left(\frac{\e_0}{\sqrt{N}C_3\|\,|\mathbf{G}|+|\mathbf{H}|\,\|_{L^q(\Theta)}}\right)^\frac{2q\omega^\prime}{q-2\omega^\prime}>0,
\end{equation}
independently of $k\in\R$. Now if $l=\sup_\Theta u>0$, we may choose $k>0$ and let $k\nearrow l$ to arrive at a contradiction. Indeed,
$v$ tends to $0$ almost everywhere in $\Theta$ as $k\nearrow l$ and also in $L^{2\omega}(\Theta)$.  By \eqref{25} we obtain that
$(v,\nabla v)$ tends to $(0,\mathbf{0})$ in $QH^1(\Theta)$. But then \eqref{26} gives
$$0=\lim_{k\rightarrow l^-}\int_\Theta\crochet{\nabla v}{Q(x)\nabla v}\,dx=
\lim_{k\rightarrow l^-}\int_{\{u>k\}}\crochet{\nabla u}{Q(x)\nabla
u}\,dx>0$$ and we have that $$\sup_\Theta u\leq0.$$ Repeating
the above argument, this time with $(v,\nabla
v)=\big((u+k)^-,\chi_{\{u<-k\}}\nabla u\big)\in QH^1(\Theta)$, see \cite{SW2},  for
any $k$ satisfying $0<k<-\inf_\Theta u$ we obtain that
$$\inf_\Theta u\geq0.$$ We conclude that $u=0$ almost everywhere on
$\Theta$.  Thus, any weak solution of problem \eqref{23} is of the
form $(u,\nabla u)=(0,\mathbf{h})\in QH^1(\Theta)$. By Definition
\ref{def2}, since $(u,\nabla u)=(0,\mathbf{h})$ is a solution of
\eqref{23}, we must have
\begin{eqnarray}
\nonumber
0\,=\,\mathcal{L}(u,u)\!&\!=\!&\!\int_\Theta\crochet{\nabla
u}{P(x)\nabla u}\,dx+\int_\Theta u\mathbf{HR}u\,dx+\int_\Theta
u\mathbf{GS}u\,dx+\int_\Theta Fu^2\,dx\\
\nonumber\!&\!\geq\!&\!
c_1\int_\Theta\crochet{\mathbf{h}}{Q(x)\mathbf{h}}\,dx.
\end{eqnarray}
Thus, $$\|u\|_{QH^1(\Theta)}=\left(\int_\Theta
u^2\,dx+\int_\Theta\crochet{\nabla u}{Q(x)\nabla
u}\,dx\right)^\frac{1}{2}=\left(\int_\Theta\crochet{\mathbf{h}}{Q(x)\mathbf{h}}\,dx\right)^\frac{1}{2}=0.$$
We conclude that $(u,\nabla u)=(0,\mathbf{0})$, completing the proof in case \eqref{24} holds.

If instead there exists $\e>0$ such that
for every $(u,\nabla u),(v,\nabla v)\in QH^1(\Theta)$ satisfying
$uv\geq0$ almost everywhere in $\Theta$ one has
\begin{equation*}
\int_\Theta Fuv+\mathbf{HR}(uv)\,dx\geq\e\int_\Theta uv\,dx,
\end{equation*}
see Definition \ref{neg1}, then we can repeat the above argument for
the adjoint homogeneous $X$--Neumann Problem
\begin{equation}\label{33}
\begin{cases}
 X^*v=0\qquad\qquad\qquad\qquad\quad\,\,\text{in }\Theta\\
 \crochet{\nu}{P(x)\nabla v+v\mathbf{G}\mathbf{S}}=0\,\qquad\text{on }\partial\Theta
\end{cases}
\end{equation}
and conclude that it admits only the trivial weak solution
$(v,\nabla v)=(0,\mathbf{0})\in QH^1(\Theta)$. Now suppose that
problem \eqref{23} has a nontrivial weak solution $(u,\nabla u)\in
QH^1(\Theta)$, then by Theorem \ref{thm2} also problem \eqref{33}
must admit a nontrivial weak solution, a contradiction. Thus the
only weak solution of problem \eqref{23} is $(u,\nabla
u)=(0,\mathbf{0})$, and the proof is complete.
\begin{flushright}$\Box$\end{flushright}

\begin{cor}\label{cor2}
Assume the hypotheses of Theorem \ref{thm5} hold.  Then the set
$\Sigma$ of real eigenvalues of the $X$-Neumann Problem
\eqref{NPeigenhom} satisfies $\Sigma\subset(0,\infty)$.
\end{cor}

\begin{cor}\label{cor3}
Assume the hypotheses of Theorem \ref{thm5} hold. Then
$0\notin\Sigma$ and the $X$-Neumann Problem \eqref{NP} admits a
unique weak solution $u\in QH^1(\Theta)$ for every $f\in
L^2(\Theta)$, every $K\in\N$, every $K$--tuple $\mathbf{T}$ of
subunit vector fields and every
$\mathbf{g}\in\big[L^2(\Theta)\big]^K$. Moreover there exists a
constant $C>0$, independent of $u$, $f$, $K$ $\mathbf{T}$ and
$\mathbf{g}$, such that $$\|u\|_{QH^1(\Theta)}\leq
C\big(\|f\|_{L^2(\Theta)}+\sqrt{K}\|\,|\mathbf{g}|\,\|_{L^2(\Theta)}\big)$$
when $u\in QH^1(\Theta)$ is a weak solution of \eqref{NP}.
\end{cor}

\textbf{Proof of Corollaries \ref{cor2} and \ref{cor3}:} These
corollaries are simple consequences of Theorems \ref{thm2},
\ref{thm3} and \ref{thm5}. \begin{flushright}$\Box$\end{flushright}

\begin{rem} Let $(\Omega,\rho)$ be a geometric homogeneous space and let $\Theta$ be a bounded domain with $\overline\Theta\subset\Omega$.  Assume that the Global Sobolev Inequality with gain $\sigma>1$ holds, see \eqref{GWP}, and that the local Poincar\'e inequality \eqref{Poincare-loc} with $p=2$ holds. Let $X$ be a second order linear degenerate
elliptic operator with rough coefficients as in \eqref{diffop} and
assume $F\in L^t(\Theta)$ with $t>\sigma^\prime$ and
$\mathbf{G},\,\mathbf{H}\in \big[L^q(\Theta)\big]^N$ with
$q>2\sigma^\prime$. It was shown in \cite{ScottCJM} that negativity condition $(2)$--i) for the operator $X$, see Definition
\ref{neg2}, is sufficient for the well--posedness of the Dirichlet
Problem
\begin{equation*}
\begin{cases}
 Xu=f+\mathbf{T}^\prime\mathbf{g}\,\qquad\text{in }\Theta\\
 u=0\qquad\qquad\qquad\text{on }\partial\Theta,
\end{cases}
\end{equation*}
with $f\in L^2(\Theta)$, $K\in\N$, $\mathbf{T}$ a $K$--tuple of
subunit vector fields and $\mathbf{g}\in\big[L^2(\Theta)\big]^K$. As
is shown by Example \ref{exe1} below, negativity condition $(2)$
for $X$ is not sufficient for the well--posedness of the
corresponding Neumann Problem \eqref{NP}.
\end{rem}

\begin{rem}
All the preceding results easily extend to include complex valued
weak solutions and complex eigenvalues/eigenfunctions.
\end{rem}

\begin{thm}\label{thm6}
Let $(\Omega,\rho)$ be a geometric homogeneous space and $\Theta$ be a bounded domain such that
$\overline\Theta\subset\Omega$.  Assume that the local Poincar\'e inequality \eqref{Poincare-loc} with $p=2$ holds and that the global weak
Poincar\'{e} inequality \eqref{globsob} with gain $\omega>1$
hold.  Let $X$ be a second order linear degenerate elliptic operator with
rough coefficients as in \eqref{diffop}. Assume that $F\in
L^t(\Theta)$ with $t>\omega^\prime$ and
that $\mathbf{G},\,\mathbf{H}\in \big[L^q(\Theta)\big]^N$ with
$q>2\omega^\prime$. If the operator $X$ is self-adjoint (see Definition \ref{def3}), then
\begin{itemize}
\item[1)] All the eigenvalues of the $X$--Neumann Problem \eqref{NPeigenhom} are real, infinite and can be ordered in a monotone sequence
which diverges to $+\infty$.
\item[2)] One has $$\lambda_1\,\,=\,\,\min\Sigma\,\,=\,\,\min_{u\in
QH^1(\Theta)\setminus\{(0,{\bf
h})\}}\frac{\mathcal{L}(u,u)}{\left(\int_\Theta u^2\,dx\right)}.$$
Moreover there exists an eigenfunction $(u_1,\nabla u_1)\in
QH^1(\Theta)$ of the Neumann Problem \eqref{NPeigenhom} related to
the eigenvalue $\lambda_1$ such that $u_1\geq0$ a.e. in $\Theta$.
\item[3)] One has that $$\lambda_2\,\,=\,\,\min\left\{\frac{\mathcal{L}(u,u)}{\left(\int_\Theta
u^2\,dx\right)}\,\bigg|\,u\in QH^1(\Theta)\setminus\{(0,{\bf
h})\},\int_\Theta uu_1\,dx=0\right\}$$ is an eigenvalue of the
Neumann Problem \eqref{NPeigenhom}, with corresponding eigenfunction
$(u_2,\nabla u_2)\in QH^1(\Theta)$ whose first component $u_2$ is
orthogonal to $u_1$ in $L^2(\Theta)$. Recursively, for every
$k\in\N$
$$\,\,\,\,\lambda_k\,\,=\,\,\min\left\{\frac{\mathcal{L}(u,u)}{\left(\int_\Theta
u^2\,dx\right)}\,\bigg|\,u\in QH^1(\Theta)\setminus\{(0,{\bf
h})\},\int_\Theta uu_j\,dx=0\,\,\text{for
all}\,j=1,\ldots,k-1\right\}$$ is an eigenvalue of the Neumann
Problem \eqref{NPeigenhom}, with corresponding eigenfunction
$(u_k,\nabla u_k)\in QH^1(\Theta)$ whose first component $u_k$ is
orthogonal to $u_j$ in $L^2(\Theta)$ for every $j=1,\ldots,k-1$.
\item[4)] $\lambda\in\R$ is
an eigenvalue if and only if $\lambda=\lambda_k$ for some $k\in\N$.
\item[5)] The sequence $\{u_k\}_{k\in\N}\subset L^2(\Theta)$ forms a
complete orthogonal system of $L^2(\Theta)$. The sequence
$\{(u_k,\nabla u_k)\}_{k\in\N}\subset QH^1(\Theta)$ is a an
independent system of elements of $QH^1(\Theta)$, which is also a
system of generators of the whole space if and only if the
projection map $i:QH^1(\Theta)\rightarrow L^2(\Theta)$ is injective.
\item[6)] The $X$-Neumann Problem \eqref{NPeigen} is variational,
with associated functional defined on $QH^1(\Theta)$ by
$$I(u)\,=\,\frac{1}{2}\mathcal{L}(u,u)-\frac{\lambda}{2}\int_{\Theta}u^2\,dx-\int_\Theta fu+\mathbf{gT}u\,dx.$$
\end{itemize}
\end{thm}

\textbf{Proof of Theorem \ref{thm6}:} The proof of this Theorem is a
standard application of functional analysis techniques, see for
instance chapter 8.12 of \cite{GT}. \begin{flushright}$\Box$\end{flushright}

\begin{exe}\label{exe1}
Let $(\Omega,\rho)$ be a geometric homogeneous space and fix a
bounded piecewise $C^1$ domain $\Theta$ with
$\overline{\Theta}\subset\Omega$. Let $P=P(x)$ be an $n\times n$
matrix as in \eqref{diffop} and consider the following Neumann
Problem
\begin{equation}\label{NPR}
\begin{cases}
 -\text{div}\big(P(x)\nabla u\big)=f+\mathbf{T}^\prime\mathbf{g}\,\,\,\qquad\text{in }\Theta\\
 \crochet{\nu}{P(x)\nabla u}=0\qquad\qquad\,\,\,\,\qquad\text{on }\partial\Theta,
\end{cases}
\end{equation}
where $f\in L^2(\Theta)$, $K\in\N$, $\mathbf{T}$ is a $K$--tuple of
subunit vector fields and $\mathbf{g}\in\big[L^2(\Theta)\big]^K$.
Assume that the following hold:
\begin{itemize}
\item[1)] the global weak Poincar\'{e} inequality on $\Theta$ with gain $\omega>1$, see
\eqref{globsob},
\item[2)] the local Poincar\'e inequality \eqref{Poincare-loc} with $p=2$.
\end{itemize}
 Theorems \ref{thm2}, \ref{thm3} and \ref{thm4} apply to
this problem where $$Xu=-\text{div}\big(P(x)\nabla
u\big)\qquad\text{ and }\qquad \mathcal{L}(u,v)=\int_\Theta
\crochet{\nabla v}{P(x)\nabla u}\,dx.$$ Moreover, the problem is self
adjoint so that Theorem \ref{thm6} also applies. Let $N\subset
QH^1(\Theta)$ be the subspace of weak solutions of
\begin{equation*}
\begin{cases}
-\text{div}\big(P(x)\nabla u\big)=0,\quad\text{in }\Theta\\
 \crochet{\nu}{P(x)\nabla u}=0,\,\,\,\,\,\,\,\quad\text{on }\partial\Theta.
\end{cases}
\end{equation*}
Obviously, $$\big\{(u,\nabla u)=(c,\mathbf{0})\in
QH^1(\Theta)\,\big|\,c\in\R\big\}\subseteq N.$$ Now suppose
$(u,\nabla u)\in N$. Then by Definition \ref{def2}, for every $v\in
QH^1(\Theta)$ one has
$$\mathcal{L}(u,v)=\int_\Theta\crochet{\nabla v}{P(x)\nabla u}\,dx=0.$$
In particular, choosing $(v,\nabla v)=(u,\nabla u)$,
$$\int_\Theta\crochet{\nabla u}{Q(x)\nabla u}\,dx\leq\frac{1}{c_1}\int_\Theta\crochet{\nabla u}{P(x)\nabla u}\,dx=0,$$
so that $(u,\nabla u)=(u,\mathbf{0})$ in $QH^1(\Theta)$. Applying the local Poincar\'{e} inequality \eqref{Poincare-loc} with $p=2$, for
every quasimetric ball $B_r(y)\subset\Theta$ with $\mathfrak{b}r\in\big(0,r_1(y)\big)$ one has
$$\int_{B_r(y)}|u-u_{B_r(y)}|\,dx\leq C\left(\int_{B_{\mathfrak{b}r}(y)}|\sqrt{Q(x)}\nabla u|^2\,dx\right)^\frac{1}{2}=0,$$
where $C>0$ is independent of $(u,\nabla u)\in QH^1(\Theta)$. Hence
$u=u_{B_r(y)}$ almost everywhere on $B_r(y)$. As quasimetric balls
are open sets, the function $u\in L^2(\Theta)$ is locally constant
in $\Theta$.  Since $\Theta$ is connected, $u$ must therefore be
constant in $\Theta$. This proves that $$\big\{(u,\nabla
u)=(c,\mathbf{0})\in QH^1(\Theta)\,\big|\,c\in\R\big\}= N.$$ By
Theorem \ref{thm2}, problem \eqref{NPR} admits a weak solution $u\in
QH^1(\Theta)$ if and only if $$\int_\Theta
cf+\mathbf{gT0}\,dx=0\qquad\text{for every }(c,\mathbf{0})\in N,$$
i.e. if and only if $$\int_\Theta f\,dx=0.$$ Now it is clear that if
$(u,\nabla u)\in QH^1(\Theta)$ is a weak solution of problem
\eqref{NPR}, then so is $(u+c,\nabla u)$ for every $c\in\R$. Hence
it is also clear that for every $f\in L^2(\Theta)$ satisfying
$\int_\Theta f\,dx=0$, every $K\in\N$, every $K$--tuple $\mathbf{T}$
of subunit vector fields and every
$\mathbf{g}\in\big[L^2(\Theta)\big]^K$ there exists a unique weak
solution $(u,\nabla u)\in QH^1(\Theta)$ of problem \eqref{NPR}
satisfying $$\int_\Theta u\,dx=0.$$

{\bf Claim:} We claim that there exists a positive constant $C$ so that
if $f\in L^2(\Theta)$ satisfies $\int_\Theta f\,dx=0$, $K\in\N$,
$\mathbf{T}$ is a $K$--tuple of subunit vector fields, $\mathbf{g}\in\big[L^2(\Theta)\big]^K$
and $(u,\nabla u)\in QH^1(\Theta)$ is a weak solution of problem
\eqref{NPR} with $\int_\Theta u\,dx=0$, then
\begin{equation}\label{32}
\|u\|_{QH^1(\Theta)}\leq
C\big(\|f\|_{L^2(\Theta)}+\sqrt{K}\|\,|\mathbf{g}|\,\|_{L^2(\Theta)}\big).
\end{equation}

To see this, notice that Proposition \ref{compact2} applies under
our current assumptions and we have that the projection onto the
first component is a compact mapping from $QH^1(\Theta)$ into
$L^2(\Theta)$.  Therefore, we are able to apply Theorem \ref{thmP}
with $p=2$, see Section 5, to conclude inequality \eqref{quasi-GP}
with $q=2$. Thus, $$\int_{\Theta}
u^2\,dx=\int_\Theta|u-u_\Theta|^2\,dx\leq
C_5\int_\Theta\crochet{\nabla u}{Q(x)\nabla u}\,dx.$$ It is now easy
to see that
\begin{equation}\label{30}
\|u\|_{QH^1(\Theta)}^2\leq(C_5+1)\int_\Theta\crochet{\nabla u}{Q(x)\nabla u}\,dx.
\end{equation}
Using that $u\in QH^1(\Theta)$ is a weak solution of \eqref{NPR} together with the subuniticity of ${\bf T}$, we have
\begin{eqnarray}
\nonumber \int_\Theta\crochet{\nabla u}{Q(x)\nabla
u}\,dx\!&\!\leq\!&\!\frac{1}{c_1}\int_\Theta\crochet{\nabla
u}{P(x)\nabla u}\,dx\\
\nonumber\!&\!=\!&\!\frac{1}{c_1}\int_\Theta fu+\mathbf{gT}u\,dx\\
\label{31}\!&\!\leq\!&\!\frac{1}{c_1}\left(\|f\|_{L^2(\Theta)}\|u\|_{L^2(\Theta)}
+\|\,|\mathbf{g}|\,\|_{L^2(\Theta)}\|\,|\mathbf{T}u|\,\|_{L^2(\Theta)}\right)\\
\nonumber\!&\!\leq\!&\!\frac{1}{c_1}\big(\|f\|_{L^2(\Theta)}+\sqrt{K}\|\,|\mathbf{g}|\,\|_{L^2(\Theta)}\big)\|u\|_{QH^1(\Theta)}.
\end{eqnarray}
Thus \eqref{30} and \eqref{31} together yield \eqref{32}, proving
our claim with $C=\frac{C_5+1}{c_1}$. Furthermore, \eqref{30} also
shows that
$$\int_\Theta\crochet{\nabla u}{Q(x)\nabla u}\,dx\leq\|u\|^2_{QH^1(\Theta)}\leq(C_5+1)\int_\Theta\crochet{\nabla u}{Q(x)\nabla
u}\,dx$$ for every $u\in QH^1(\Theta)$ with $\int_\Theta u\,dx=0$.
Therefore, the norm
$$\|u\|_{QH^1_*(\Theta)}=\int_\Theta\crochet{\nabla u}{Q(x)\nabla
u}\,dx,$$ is equivalent to $\|\cdot\|_{QH^1(\Theta)}$ on the
subspace $QH^1_*(\Theta)\subset QH^1(\Theta)$ defined by
$$QH^1_*(\Theta)\,=\,\left\{u\in QH^1(\Theta)\,\bigg|\,\int_\Theta
u\,dx=0\right\}\,=\,N^\perp.$$


We also mention that Theorem \ref{thm5} shows that the problem
\begin{equation}\label{NPReigen}
\begin{cases}
-\text{div}\big(P(x)\nabla u\big)=\lambda u+f+\mathbf{T}^\prime\mathbf{g}\,\,\qquad\text{in }\Theta\\
 \crochet{\nu}{P(x)\nabla u}=0\qquad\,\qquad\qquad\qquad\,\,\,\,\text{on }\partial\Theta,
\end{cases}
\end{equation}
with $\lambda<0$ admits a unique weak solution $u\in QH^1(\Theta)$
for every $f\in L^2(\Theta)$, every $K\in\N$,  every $K$--tuple
$\mathbf{T}$ of subunit vector fields and every
$\mathbf{g}\in\big[L^2(\Theta)\big]^K$. Hence, all eigenvalues of
problem \eqref{NPReigen} must be nonnegative with
$\lambda_1=\min\Sigma=0$ where $\Sigma\subset\R$ is the set of
eigenvalues of problem \eqref{NPReigen}.  Furthermore, the
eigenvalue $\lambda_1$ is simple. Since problem \eqref{NPReigen} is
self-adjoint, its eigenvalues form a monotone sequence diverging to
$+\infty$. The corresponding eigenfunctions $\{(u_k,\nabla
u_k)\}_{k\in\N}\subset QH^1(\Theta)$ form an independent system of
elements such that $\{u_k\}_{k\in\N}\subset L^2(\Theta)$ is a
complete orthogonal system. Moreover, one can choose $(u_1,\nabla
u_1)=(1,\mathbf{0})\in QH^1(\Theta)$. If the projection map
$i:QH^1(\Theta)\rightarrow L^2(\Theta)$ is injective, then the
eigenfunctions $\{(u_k,\nabla u_k)\}_{k\in\N}$ are also a system of
generators of $QH^1(\Theta)$.

One has that $(u,\nabla u)\in QH^1(\Theta)$ is a solution of problem
\eqref{NPReigen} if and only if it is a critical point of the
functional $I:QH^1(\Theta)\rightarrow\R$ defined by
$$I(u)\,=\,\frac{1}{2}\int_\Theta\crochet{\nabla u}{P(x)\nabla u}\,dx-\frac{\lambda}{2}\int_\Theta u^2\,dx-\int_\Theta fu+\mathbf{gT}u\,dx.$$

Finally, if $P(x)=Q(x)$ almost everywhere in $\Theta$ one can show
that $\{(u_k,\nabla u_k)\}_{k\in\N}\subset QH^1(\Theta)$ is an
orthogonal system. Indeed, for every $k\neq j$,
$k,j\in\N$, one has
$$\int_\Theta\crochet{\nabla u_k}{Q(x)\nabla u_j}\,dx=\mathcal{L}(u_j,u_k)=\lambda_j\int_\Theta
u_ju_k\,dx=0.$$ Hence for every $k\neq j$, $k,j\in\N$, we obtain
$$(u_j,u_k)_{QH^1(\Theta)}=\int_\Theta\crochet{\nabla u_j}{Q(x)\nabla u_k}\,dx+\int_\Theta u_ju_k\,dx=0,$$
as claimed.

The results on the Neumann Problem \eqref{NPR} described above extend and complete previous results on such problems obtained in \cite{ScottPhD}.

\end{exe}
\vspace{0,3cm}

\section{Spectral Results for the $X$-Dirichlet Problem}

In this section we give a spectral theorem related to $X$-Dirichlet
problems (that is, Dirichlet problems associated to the operator $X$
with rough coefficients as in (\ref{diffop}). Existence results for $X$-Dirichlet problems are given in
\cite{ScottCJM} and we refer the interested reader there for
statements and proofs.  We begin by recalling the definition of weak
solution related to the $X$-Dirichlet Problem on a bounded domain
with homogeneous boundary data as given in \cite{ScottCJM}.

\begin{defn}\label{def4}
Let $X$ be a second order linear degenerate elliptic operator with
rough coefficients as in \eqref{diffop}. Assume that the global
Sobolev inequality with gain $\sigma>1$ holds, see \eqref{GWP}.  Let
$\mathbf{G},\,\mathbf{H}\in\big[L^q(\Theta)\big]^N$ with
$q\geq2\sigma^\prime$ and let $F\in L^t(\Theta)$ with
$t\geq\sigma^\prime$. If $f\in L^2(\Theta)$, $K\in\N$, $\mathbf{T}$
is a $K$--tuple of subunit vector fields and
$\mathbf{g}\in\big[L^2(\Theta)\big]^K$, a function $(u,\nabla u)\in
QH^1_0(\Theta)$ is a weak solution of the Dirichlet Problem
\begin{equation}\label{DP}
\begin{cases}
 Xu=f+\mathbf{T^\prime}\mathbf{g}\,\,\quad\text{in }\Theta\\
 u=0\,\,\,\,\,\,\,\quad\quad\qquad\text{on }\partial\Theta
\end{cases}
\end{equation}
if and only if
\begin{equation}\label{27}
\mathcal{L}(u,v)=\int_\Theta
fv+\mathbf{g}\mathbf{T}v\,dx\;\text{ for all}\,\; v\in QH_0^1(\Theta).
\end{equation}
\end{defn}

\begin{thm}\label{thm7}
Let $(\Omega,\rho)$ be a geometric homogeneous space and let $\Theta$ be a bounded domain such that
$\overline\Theta\subset\Omega$.  Assume that the local Poincar\'e inequality \eqref{Poincare-loc} holds with $p=2$ and that the global Sobolev
inequality \eqref{GWP} with gain $\sigma>1$ holds. Let $X$ be a second order linear degenerate
elliptic operator with rough coefficients as in
\eqref{diffop}. Assume that $F\in L^t(\Theta)$ with
$t>\sigma^\prime$ and that
$\mathbf{G},\,\mathbf{H}\in \big[L^q(\Theta)\big]^N$ with
$q>2\sigma^\prime$. Then each of the following hold.
\begin{itemize}
\item[1)] There exists an at most countable set $\Sigma\subset\R$
such that the $X$-Dirichlet problem
\begin{equation}\label{DPeigen}
\begin{cases}
 Xu=\lambda u+f+\mathbf{T}^\prime\mathbf{g}\,\quad\text{in }\Theta\\
 u=0\qquad\qquad\qquad\quad\,\text{on }\partial\Theta
\end{cases}
\end{equation}
admits a unique weak solution $u\in QH^1_0(\Theta)$ for every $f\in
L^2(\Theta)$, every $K\in\N$, every $K$--tuple $\mathbf{T}$ of
subunit vector fields and every
$\mathbf{g}\in\big[L^2(\Theta)\big]^K$ if and only if
$\lambda\notin\Sigma$.
\item[2)] If $\Sigma$ is infinite, its elements can be arranged in a
monotone sequence diverging to $+\infty$.
\item[3)] If $\lambda\notin\Sigma$ there exists a constant
$C=C(\lambda,\Theta,\Omega,c_1,C_1,\mathbf{G},\mathbf{H},F)>0$ such
that $$\|u\|_{QH^1_0(\Theta)}\leq
C\big(\|f\|_{L^2(\Theta)}+\sqrt{K}\|\,|\mathbf{g}|\,\|_{L^2(\Theta)}\big)$$
whenever $f\in L^2(\Theta)$, $K\in\N$, $\mathbf{T}$ is a $K$--tuple
of subunit vector fields, $\mathbf{g}\in\big[L^2(\Theta)\big]^K$ and
$u\in QH^1_0(\Theta)$ is a weak solution of \eqref{DP}.
\item[4)] If $\lambda\in\Sigma$, let $N\subset QH^1_0(\Theta)$ be
the subspace of weak solutions of the $X$-Dirichlet Problem
\begin{equation*}
\begin{cases}
 Xu=\lambda u\quad\text{in }\Theta\\
 u=0\qquad\,\,\,\text{on }\partial\Theta,
\end{cases}
\end{equation*}
and let $N^*\subset QH^1_0(\Theta)$ be the subspace of weak
solutions of the adjoint problem
\begin{equation*}
\begin{cases}
 X^*u=\lambda u\qquad\text{in }\Theta\\
 u=0\qquad\qquad\text{on }\partial\Theta.
\end{cases}
\end{equation*}
Then $1\leq\text{dim}\,N=\text{dim}\,N^*<\infty$ and problem
\eqref{DP} admits a weak solution $u\in QH^1_0(\Theta)$ if and only
if
$$\int_\Theta fv+\mathbf{gT}v\,dx=0\;\text{ for all}\;\,v\in N^*.$$
\item[5)] If $X$ satisfies negativity condition $(2)$, see Definition \ref{neg2}, then
$\Sigma\subset(0,\infty)$.
\item[6)] If $X$ is self--adjoint (that is, if $\mathbf{HR}=\mathbf{GS}$ almost everywhere in
$\Theta$), then all eigenvalues of $X$ are real, $\Sigma$ is
infinite and we have the following variational characterization of
the eigenvalues of $X$:
$$\lambda_1\,\,=\,\,\min\Sigma\,\,=\,\,\min_{u\in
QH^1_0(\Theta)\setminus\{(0,{\bf
h})\}}\frac{\mathcal{L}(u,u)}{\left(\int_\Theta u^2\,dx\right)},$$
and there exists an eigenfunction $(u_1,\nabla u_1)\in
QH^1_0(\Theta)$ of the $X$-Dirichlet Problem \eqref{DP} related to
the eigenvalue $\lambda_1$ for whom $u_1\geq0$ a.e. in $\Theta$.
Furthermore,
$$\lambda_2\,\,=\,\,\min\left\{\frac{\mathcal{L}(u,u)}{\left(\int_\Theta
u^2\,dx\right)}\,\bigg|\,u\in QH^1_0(\Theta)\setminus\{(0,{\bf
h})\},\int_\Theta uu_1\,dx=0\right\},$$  with corresponding
eigenfunction $(u_2,\nabla u_2)\in QH^1_0(\Theta)$ where $u_2$ is
orthogonal to $u_1$ in $L^2(\Theta)$.  Recursively, for every
$k\in\N$
$$\qquad\,\,\lambda_k\,\,=\,\,\min\left\{\frac{\mathcal{L}(u,u)}{\left(\int_\Theta
u^2\,dx\right)}\,\bigg|\,u\in QH^1_0(\Theta)\setminus\{(0,{\bf
h})\},\int_\Theta uu_j\,dx=0\,\,\text{for
all}\,j=1,\ldots,k-1\right\},$$  with corresponding eigenfunction
$(u_k,\nabla u_k)\in QH^1_0(\Theta)$, where $u_k$ is orthogonal to
$u_j$ in $L^2(\Theta)$ for every $j=1,\ldots,k-1$. Moreover,
$\lambda\in\R$ is an eigenvalue if and only if $\lambda=\lambda_k$
for some $k\in\N$. The sequence $\{u_k\}_{k\in\N}\subset
L^2(\Theta)$ forms a complete orthogonal system of $L^2(\Theta)$.
The sequence $\{(u_k,\nabla u_k)\}_{k\in\N}\subset QH^1_0(\Theta)$
is an independent system of elements of $QH^1_0(\Theta)$, which is
also a system of generators of $QH^1_0(\Theta)$ if and only if the
projection map $i:QH^1_0(\Theta)\rightarrow L^2(\Theta)$ is
injective. Finally, problem \eqref{DPeigen} is variational with
associated functional defined on $QH^1_0(\Theta)$ by
$$I(u)=\frac{1}{2}\mathcal{L}(u,u)-\frac{\lambda}{2}\int_\Theta u^2\,dx-\int_\Theta fu+\mathbf{gT}u\,dx.$$
\end{itemize}
\end{thm}

\textbf{Proof of Theorem \ref{thm7}:} The proof of this theorem is
similar to the proofs of the preceding section where the global Poincar\'e inequality \eqref{globsob} with gain $\omega>1$ is replaced with the global Sobolev inequality \eqref{GWP} when necessary.  We therefore omit the proof. \begin{flushright}$\Box$\end{flushright}

\section{A Maximum Principle for Second Order Linear Degenerate Elliptic Equations with Rough Coefficients}

This section contains a maximum principle for weak solutions of the
differential inequality $Xu\leq 0$ for second order degenerate
elliptic operators $X$ with rough coefficients.  To this end, we fix
a geometric homogeneous space $(\Omega,\rho)$ with $\Omega$ as in
Section 1.  We also fix an $n\times n$ matrix $Q(x)$ as in Section 1
and let $\Theta$ be a bounded domain satisfying
$\overline{\Theta}\subset\Omega$.  Furthermore, in order to avoid
confusion we will refer to an element of $QH^1(\Theta)$ by writing
${\bf u}\in QH^1(\Theta)$ where ${\bf u}=(u,\nabla u)$.  We begin by
giving the definition of weak solution of $X\bu\leq 0$ in $\Theta$.

\begin{defn}\label{def5}
We say that $\bu\in QH^1(\Theta)$ is a weak solution of
\bea\label{ineq} X\bu\leq 0\quad\text{in }\Theta
\eea
if and only if
$\mathcal{L}(\bu,\bv)\leq0$ for every $\bv=(v,\nabla v)\in QH^1_0(\Theta)$
satisfying $v\geq0$ almost everywhere in $\Theta$.
\end{defn}

In order to state our maximum principle, we define a notion of
non-positivity for the first component $u$ of an element
$\bu=(u,\nabla u)\in QH^1(\Theta)$ in terms of membership in the
space $QH^1_0(\Theta)$.  At this time it is useful to recall that if
$\bu=(u,\nabla u)\in QH^1(\Theta)\;\big(QH^1_0(\Theta)\text{
respectively}\big)$ then, as is shown in \cite{SW2}, $\bu^+ =
(u^+,\chi_{\{u>0\}}\nabla u)\in QH^1(\Theta)$
$\big(QH^1_0(\Theta)\text{ respectively}\big)$.

\begin{defn}\label{def6}\hspace{1in}\begin{enumerate}
\item We say that $\bu=(u,\nabla u)\in QH^1(\Theta)$ satisfies $u\leq0$ on
$\partial\Theta$ if and only if $$\bu^+=(u^+,\chi_{\{u>0\}}\nabla u)\in QH^1_0(\Theta).$$
\item We say that ${\bf u},{\bf v}\in QH^1(\Theta)$ satisfy $ u\leq
 v$ on $\partial\Theta$ if and only if $u-v\leq0$ on
$\partial\Theta$ in the sense of item $(1)$.
\item For each $k\in \mathbb{R}$ recall that ${\bf k}=(k,0)\in QH^1(\Theta)$ as $\Theta$ is bounded.  Thus, for $\bu=(u,\nabla u)\in QH^1(\Theta)$ we define
$$\sup_{\partial\Theta}{ u}=\inf\big\{k\in\R\,\,\big|\,\, u\leq k\text{ on }\partial\Theta\big\},\;\text{and}\;
\inf_{\partial\Theta}u=-\sup_{\partial\Theta}(-u).$$
\end{enumerate}
\end{defn}

\begin{thm}\label{thm8}
Noting the first paragraph of this section, assume that the local Poincar\'e inequality \eqref{Poincare-loc}
with $p=2$ holds and that the global Sobolev inequality \eqref{GWP}
with gain $\sigma>1$ holds. Let $X$ be a second order linear
degenerate elliptic operator with rough coefficients as in
\eqref{diffop} that satisfies negativity condition (2)--i), see
Definition \ref{neg2}. Assume that $F\in L^t(\Theta)$ with
$t>\sigma^\prime$ and that $\mathbf{G},\,\mathbf{H}\in
\big[L^q(\Theta)\big]^N$ with $q>2\sigma^\prime$. If ${\bf u}\in
QH^1(\Theta)$ is a weak solution of \eqref{ineq} then
$$\sup_\Theta u\leq\sup_{\partial\Theta}u^+.$$
\end{thm}

\textbf{Proof of Theorem \ref{thm8}:} We argue by contradiction. Let
$\displaystyle l=\sup_{\Theta}u$, let $\displaystyle
m=\sup_{\partial\Theta}u^+$ and suppose that $l>m$. Fix
$k\in\mathbb{R}$ with $m< k<l$ and set ${\bv_k}=(v_k,\nabla
v_k)=\big((u-k)^+,\chi_{\{u>k\}}\nabla u\big)\in QH^1(\Theta)$.
Since $k>m$ and since $m\geq0$, it is not difficult to see that
$\bv_k\in QH^1_0(\Theta)$ by Definition \ref{def6}.  Moreover,
$v_k\geq0$ almost everywhere in $\Theta$. Thus, $uv_k\geq0$ almost
everywhere in $\Theta$, and an application of Lemma \ref{lem1}
implies that
$\mathbf{GS}({\bu\bv_k})=u\mathbf{GS}{\bv_k}+v_k\mathbf{GS}{\bu}$
and $\mathbf{HR}({\bu\bv_k})=u\mathbf{HR}{\bv_k}+v_k\mathbf{HR}{\bf
u}$.

Consider now the case where $\mathbf{G}=\mathbf{H}=0$ almost
everywhere on $\Theta$. Using that $X$ satisfies negativity
condition $(2)$--i) and that $\nabla v_k=\nabla u$ a.e. on the
support of $\nabla v_k$, the definition of weak solution of
\eqref{ineq} yields
\begin{eqnarray}
\nonumber0\,\,\geq\,\,\mathcal{L}({\bu},{\bv_k})\!&\!=\!&\!\int_{\Theta}\crochet{\nabla
v_k}{P(x)\nabla u}\,dx+\int_\Theta
Fuv_k\,dx\\
\nonumber\!&\!\geq\!&\!\int_{\Theta}\crochet{\nabla v_k}{P(x)\nabla
v_k}\,dx\,\,\geq\,\, c_1\int_\Theta\crochet{\nabla v_k}{Q(x)\nabla
v_k}\,dx.
\end{eqnarray}
By the global Sobolev inequality \eqref{GWP} we see that
${\bv_k}=(0,\mathbf{0})$ in $QH^1_0(\Theta)$. This in turn implies
that $v_k=(u-k)^+=0$ in $L^2(\Theta)$ and hence that $u\leq k$
almost everywhere in $\Theta$, contradicting our assumption
$k<\sup_\Theta u$.

We now focus on the case where
$\|\,|\mathbf{G}|+|\mathbf{H}|\,\|_{L^q(\Theta)}\neq0$. Since $X$
satisfies negativity condition $(2)$--i), we argue as in the proof
of Theorem \ref{thm5} to obtain
\begin{eqnarray}
  \nonumber\int_\Theta\crochet{\nabla v_k}{P(x)\nabla v_k}\,dx&=&\int_\Theta\crochet{\nabla v_k}{P(x)\nabla u}\,dx\\
  \nonumber&\leq&\int_\Theta|v_k|\big(|\mathbf{G}||\mathbf{S}{\bv_k}|+|\mathbf{H}||\mathbf{R}{\bv_k}|\big)\,dx.
\end{eqnarray}
Recalling the result of \cite[Lemma 3.18]{ScottCJM} and setting
$\Gamma=\text{supp}\big(|\sqrt{Q(x)}\nabla v_k|\big)$, we see that
\begin{eqnarray}
 \nonumber\|{\bv_k}\|^2_{QH_0^1(\Theta)}\!&\!\leq\!&\!\frac{1}{c_1}\int_\Theta\crochet{\nabla v_k}{P(x)\nabla v_k}\,dx\\
 \nonumber&\leq&\frac{\sqrt{N}}{c_1}\big\|\big(|\mathbf{G}|+|\mathbf{H}|\big){ v_k}\big\|_{L^2(\Gamma)}\left(\int_\Theta\crochet{\nabla v_k}{Q(x)\nabla
 v_k}\,dx\right)^\frac{1}{2}\\
 \nonumber&\leq&\frac{\sqrt{N}}{c_1}\|\,|\mathbf{G}|+|\mathbf{H}|\,\|_{L^q(\Theta)}\|v_k\|_{L^\frac{2q}{q-2}(\Gamma)}\|{\bv_k}\|_{QH_0^1(\Theta)}\\
 \label{28}&\leq&\frac{\sqrt{N}}{c_1}\|\,|\mathbf{G}|+|\mathbf{H}|\,\|_{L^q(\Theta)}|\Gamma|^{\frac{1}{2\sigma^\prime}-\frac{1}{q}}
 \|v_k\|_{L^{2\sigma}(\Theta)}\|{\bv_k}\|_{QH_0^1(\Theta)}\\
 \nonumber&\leq&\frac{\sqrt{N}\wt{C}}{c_1}\|\,|\mathbf{G}|+|\mathbf{H}|\,\|_{L^q(\Theta)}|\Gamma|^{\frac{1}{2\sigma^\prime}-\frac{1}{q}}
 \|{\bv_k}\|^2_{QH_0^1(\Theta)}.
\end{eqnarray}
Dividing through by $||\bv_k||^2_{QH^1_0(\Theta)}$ we have
\begin{equation}\label{29}
|\Gamma|\geq\left(\frac{c_1}{\sqrt{N}\wt{C}\|\,|\mathbf{G}|+|\mathbf{H}|\,\|_{L^q(\Theta)}}\right)^\frac{2q\sigma^\prime}{q-2\sigma^\prime}>0
\end{equation}
independently of $k$, for $m< k<l$. As in Theorem \ref{thm5},
sending $k\rightarrow l^-$ gives a contradiction and we conclude that  
$$\sup_\Theta u\leq\sup_{\partial\Theta}u^+.$$
\begin{flushright}$\Box$\end{flushright}

\section{Poincar\'e Inequalities and Compact Projection of Sobolev Spaces}

In this section we give a result demonstrating a global Poincar\'e inequality with gain as a consequence of a compact ``embedding"-type property for degenerate Sobolev spaces.  We say that the \emph{Compact Projection Property} from $QH^{1,p}(\Theta)$ into $L^q(\Theta)$ holds if and only if the projection map $i:QH^{1,p}(\Theta)\ra L^q(\Theta)$ defined by $i((u,\nabla u)=u$ is a compact mapping.  Recall that the space $QH^{1,p}(\Theta)$ is defined as the closure in the norm $$||w||_{QH^{1,p}(\Theta)} = ||w||_{L^p(\Theta)} + ||\sqrt{Q(x)}\nabla w||_{L^p(\Theta)}$$ of the collection $Lip_{Q,p}(\Theta)$ of locally Lipschitz functions defined in $\Theta$ with finite $QH^{1,p}(\Theta)$ norm. Note that $QH^{1,p}(\Theta)$ is denoted by $W^{1,p}_Q(\Theta)$ in \cite{MRW}, \cite{CRW}, and as $W^{1,p}_{\mathcal{Q}}(\Theta)$ in \cite{SW2}.

\begin{thm}\label{thmP}
Let $(\Omega,\rho)$ be a geometric homogeneous space, $\Theta$ be a bounded domain such that
$\overline\Theta\subset\Omega$ and fix $p,q$ with $1\leq p\leq q<\infty$. Assume that the
local Poincar\'{e} inequality \eqref{Poincare-loc} of order $p$ holds. If the Compact Projection Property from $QH^{1,p}(\Theta)$ to
$L^q(\Theta)$ holds, there exists a positive constant $C_5>0$
such that
\begin{equation}\label{quasi-GP}
\left(\int_\Theta |w - w_\Theta|^r\,dx\right)^\frac{1}{r} \leq
C_5\left(\int_\Theta |\sqrt{Q}\nabla w|^p\,dx\right)^\frac{1}{p}
\end{equation}
for every $r\in [1,q]$ and pair $(w,\nabla w)\in QH^{1,p}(\Theta)$.
\end{thm}

\begin{cor} Assume that the conditions of Theorem \ref{thmP} hold with $p=2$. Then,
the global Poincar\'{e} inequality \eqref{GP} holds (and also the global weak Poincar\'e inequality (\ref{globsob})) with gain
$\omega=\frac{q}{2}$.
\end{cor}

\begin{rem}
When put in context with Remark \ref{rem1.5}, Proposition \ref{compact2} (when $p=2$) and related results of \cite{CRW}, one can see a clear relationship, vis-a-vis almost necessity and sufficiency, between the compact projection property and the Poincar\'e inequality (\ref{quasi-GP}) in the setting of a geometric homogeneous space where a local Poincar\'e inequality of the form (\ref{Poincare-loc}) holds.
\end{rem}

\textbf{Proof of Theorem \ref{thmP}:} We argue by contradiction and
begin with the case where $r=q\geq p$. If the result does not
hold, then for every $k\in\N$ there is a $w_k\in QH^{1,p}(\Theta)$ so
that
\begin{equation}\label{2}
\left(\int_\Theta |w_k - (w_k)_\Theta|^q\,dx\right)^\frac{1}{q}
>k\left(\int_\Theta |\sqrt{Q}\nabla w_k|^p\,dx\right)^\frac{1}{p},
\end{equation}
where $(w_k)_\Theta=\frac{1}{|\Theta|}\int_\Theta w_k\,dx$. Thus,
$\|w_k - (w_k)_\Theta\|_{L^q(\Theta)}>0$ and we are able to define
$$v_k=\frac{w_k - (w_k)_\Theta}{\|w_k -
(w_k)_\Theta\|_{L^q(\Theta)}}.$$ It is clear that both
$$\|v_k\|_{L^q(\Theta)}=1,\quad (v_k)_\Theta=0,$$ and $$\left(\int_\Theta |\sqrt{Q}\nabla v_k|^p\,dx\right)^\frac{1}{p}<\frac{1}{k}$$ hold.
Since $q\geq p$ and $\Theta$ is bounded, the
sequence $\{(v_k,\nabla v_k)\}_{k\in\N}$ is bounded in $QH^{1,p}(\Theta)$. By the
Compact Projection Property there exists $v\in L^q(\Theta)$
so that, up to a subsequence, $v_k\rightarrow v$ pointwise a.e. in $\Theta$ and also in $L^q(\Theta)$.  It is not difficult to see that the limit $v$ satisfies
\begin{equation}\label{3}
\|v\|_{L^q(\Theta)}=1,\quad v_\Theta=0.
\end{equation}
Applying the local Poincar\'e inequality \eqref{Poincare-loc} to $v_k$ gives
\begin{eqnarray}
\label{5-sr1}\frac{1}{|B_r|}\int_{B_r} |v_k- (v_k)_{B_r}|\, dx&\le& C_2
r \left(\frac{1}{|B_{\mathfrak{b}r}|}\int_{B_{\mathfrak{b}r}}
|\sqrt{Q}\nabla v_k|^p\, dx
\right)^{\frac{1}{p}}\,\,<\,\,\frac{C_2r}{|B_{\mathfrak{b}r}|^\frac{1}{p}}\,\frac{1}{k}
\end{eqnarray}
for every $k\in\N$ and every $B_r=B(y,r)$ such that
$\mathfrak{b}r\in(0,r_1(y))$ and
$\overline{B_{\mathfrak{b}r}}\subset\Theta$. Passing to the limit as
$k$ tends to $\infty$ we get
$$\|v-v_{B_r}\|_{L^1(B_r)}=0$$ and so $v$ is $a.e.$ constant on
$B_r$. Since quasimetric balls are open, $v$ is locally constant in the connected set $\Theta$ and we conclude that $$v\equiv const.\quad\text{a.e. in }\Theta.$$ This contradicts
\eqref{3} and proves \eqref{quasi-GP} when $r=q$. Since $\Theta$ is bounded, we can recover inequality
\eqref{quasi-GP} in the case $r\in[1,q)$ by a simple application of H\"{o}lder's inequality.
\begin{flushright}$\Box$\end{flushright}

\begin{rem} Theorem \ref{thmP} is improved by replacing the local Poincar\'e inequality of order $p\geq 1$ with its weaker $L^1\rightarrow L^p$ counterpart obtained by replacing $p$ with $1$ on the left-hand side of (\ref{Poincare-loc}).  The proof of this improved version is identical to the one just given n.b. (\ref{5-sr1}).  
\end{rem}


\end{document}